\documentclass[draftcls,onecolumn,journal]{IEEEtran}
\makeatletter
\def\endthebibliography{%
	\def\@noitemerr{\@latex@warning{Empty `thebibliography' environment}}%
	\endlist
}

\usepackage{array}
\newcolumntype{L}[1]{>{\raggedright\let\newline\\\arraybackslash\hspace{0.1pt}}m{#1}}
\newcolumntype{C}[1]{>{\centering\let\newline\\\arraybackslash\hspace{0.1pt}}m{#1}}
\newcolumntype{R}[1]{>{\raggedleft\let\newline\\\arraybackslash\hspace{0.1pt}}m{#1}}

\usepackage{microtype}
\usepackage[font=scriptsize]{caption}
\usepackage[noadjust]{cite}
\usepackage{amsmath,amssymb,amsfonts}

\usepackage{graphicx}
\usepackage{comment}
\usepackage{textcomp}
\usepackage[dvipsnames]{xcolor}
\usepackage{amsbsy}
\usepackage{amsthm}
\usepackage{algorithmicx}
\usepackage[]{algorithm}
\usepackage{algpseudocode}
\usepackage{pdfpages}
\usepackage{tikz}
\usepackage{booktabs, multirow} 
\usepackage{soul}
\usepackage{changepage,threeparttable} 
\usepackage{subcaption}
\usepackage{svg, mathtools}
\usepackage{xr}
\usepackage{xcolor}
\usepackage{hyperref}

\usepackage{pifont}
\newcommand{\cmark}{\ding{51}}%
\newcommand{\xmark}{\ding{55}}%

\def \x {{\b{x}}}
\def \v {{\b{v}}}
\def \w {{\b{w}}}

\def \y {{\b{y}}}
\def \z {{\b{z}}}
\def \g {{\b{g}}}

\def \a {{\b{a}}}
\def \q {{\b{q}}}

\def \kb {{\bar{k}}}
\def \Hc {{\mathcal H}}

\def \tx {{\tb{x}}}

\def \h {{\mathbf{h}}}

\def \ft {{\tilde{f}}}

\def \B {{\b{B}}}

\def \I {{\b{I}}}

\def \T {{\mathsf{T}}}

\def \cB {{\c{B}}}

\def \cX {{\c{X}}}

\def \sX {{\mathsf{X}}}
\def \sY {{\mathsf{Y}}}

\def \xib {{\bs{\xi}}}
\def \nub {{\bs{\nu}}}
\def \nubh {{\hat{\nub}}}
\def \nubt {{\tilde{\nub}}}
\def \pib {{\bs{\pi}}}
\def \ep {{\bs{\varepsilon}}}

\def \EE {{\mathbb{E}}}
\def \Rn {{\mathbb{R}}}

\def \O {{\mathcal{O}}}

\def \e {{\ep}}

\def \xh {{\hat{\x}}}
\def \xhat {{\xh(\x_t,\z_t,\xib_t)}}
\def \Xc {{\mathcal X}}
\def \fh {{\hat{f}}}

\def \del {{\boldsymbol{\delta}}}

\def \xh {{\hat{\x}(\x_t^{s+1},\tx^s)}}

\def \cK {{\mathcal{K}}}

\def \ft {{\tilde{f}}}
\def \fh {{\hat{f}}}
\def \gt {{\tilde{g}}}
\def \gbt {{\tilde{\g}}}
\def \cJ {{\bar{\mathcal{J}}}}
\def \cI {{\bar{\mathcal{I}}}}
\def \d {{\b{d}}}
\def \lam {{\boldsymbol{\lambda}}}
\def \lamh {{\hat{\lam}}}
\def \lamt {{\tilde{\lam}}}
\def \cD {{\c{D}}}

\def \xh {{\hat{\x}}}
\def \xhat {{\xh_t}} 
\def \Xc {{\mathcal X}}

\def \fh {{\hat{f}}}

\def \w {{\mathbf{w}}}
\def \Rn {{\mathbb{R}}}

\def \Rn {{\mathbb{R}}}

\def \O {{\mathcal{O}}}
\def \Ot {{\tilde{\mathcal{O}}}}

\def \Hc {{\mathcal{H}}}

\newtheorem{assumption}{}

\theoremstyle{remark}
\newtheorem{rem}{\bf Remark}
\providecommand{\ip}[2]{\langle #1, #2 \rangle} 
\newtheorem{theorem}{Theorem}
\newtheorem{lemma}{Lemma}

\newtheorem{corollary}{Corollary}

\providecommand{\abs}[1]{\lvert#1\rvert}
\renewcommand{\b}[1]{\ensuremath{\mathbf{#1}}} 
\renewcommand{\c}[1]{\ensuremath{\mathcal{#1}}} 
\providecommand{\ip}[1]{\langle#1\rangle} 

\providecommand{\abs}[1]{\lvert#1\rvert}
\renewcommand{\b}[1]{\ensuremath{\mathbf{#1}}} 
\newcommand{\bs}[1]{\ensuremath{\boldsymbol{#1}}} 
\renewcommand{\c}[1]{\ensuremath{\mathcal{#1}}} 
\newcommand{\E}[1]{\ensuremath{\EE\left[#1\right]}}  
\newcommand{\ind}{1\hspace{-1.6mm}1} 
\newcommand{\norm}[1]{\ensuremath{\left\|#1\right\|}} 
\newcommand{\tb}[1]{\ensuremath{\tilde{\mathbf{#1}}}} 
\newcommand{\mat}[1]{\ensuremath{\begin{bmatrix}#1\end{bmatrix}}} 
\newcommand{\eqtext}[1]{\ensuremath{\stackrel{\text{#1}}{=}}} 
\newcommand{\leqtext}[1]{\ensuremath{\stackrel{\text{#1}}{\leq}}} 
\newcommand{\geqtext}[1]{\ensuremath{\stackrel{\text{#1}}{\geq}}} 
\providecommand{\ip}[2]{\langle #1, #2 \rangle} 

\begin{document}
	
	\title{ Constrained Stochastic Recursive Momentum Successive Convex Approximation
		\author{Basil M. Idrees, \IEEEmembership{Student Member, IEEE}, Lavish Arora, \IEEEmembership{Student Member, IEEE}, and Ketan Rajawat, \IEEEmembership{Member, IEEE}
			\thanks{An abridged version of this work, specifically the statements of Lemmas~\ref{le4}--\ref{dualvari}, the proof of Theorem~\ref{th2}, and some experiments from Section~\ref{ssec:eematp}, was presented at the 2024 IEEE 34th International Workshop on Machine Learning for Signal Processing (MLSP), London, UK~\cite{idrees2024non}. All other content, including additional proofs and results, constitutes contribution of this paper/manuscript.
				
	}}}
	
	\maketitle
	
	\begin{abstract}
		We consider stochastic optimization problems with non-convex functional constraints, such as those arising in trajectory generation, sparse approximation, and robust classification. To this end, we put forth a recursive momentum-based accelerated successive convex approximation (SCA) algorithm. At each iteration, the proposed algorithm entails constructing convex surrogates of the stochastic objective and the constraint functions, and solving the resulting convex optimization problem. A recursive update rule is employed to track the gradient of the stochastic objective function, which contributes to variance reduction and hence accelerates the algorithm convergence. A key ingredient of the proof is a new parameterized version of the standard Mangasarian-Fromowitz Constraints Qualification, that allows us to bound the dual variables and hence obtain problem-dependent bounds on the rate at which the iterates approach an $\epsilon$-stationary point. Remarkably, the proposed algorithm achieves near-optimal stochastic first-order (SFO) complexity with adaptive step sizes closely matching  that achieved by state-of-the-art stochastic optimization algorithms for solving unconstrained problems. As an example, we detail an obstacle-avoiding trajectory optimization problem that can be solved using the proposed algorithm and show that its performance is superior to that of the existing algorithms used for trajectory optimization. The performance of the proposed algorithm is also shown to be comparable to that of a specialized sparse classification algorithm applied to a binary classification problem.
	\end{abstract}
	
	
	\section{Introduction}
	\subsection{Background} \label{sec:back}
	In this work, we consider the stochastic non-convex constrained optimization problem:
	\begin{align}\label{pc}
		F^\star = \min_{\x \in \Rn^n} & ~ F(\x):= U(\x) + u(\x) \tag{$\mathcal{P}$}\\
		\text{s. t. } &~~~~~~  h_i(\x) \leq 0, & i = 1, \ldots, I  \nonumber  \\
		&~~~~~~ g_j(\x) \leq 0, & j = 1, \ldots, J	\nonumber
	\end{align}
	where $U(\x) := \EE[f(\x,\xib)] $ and $f$ is a smooth possibly non-convex function, while $\{g_j\}_{j=1}^J$ are smooth but possibly non-convex functions. On the other hand, the functions $\{h_i\}_{i=1}^I$ are smooth convex functions, and $u$ is a convex function. Here, the expectation is with respect to the random variable $\xib$ with an unknown distribution. We seek to solve \eqref{pc} using stochastic approximation schemes that require independent samples of $\xib$ observed sequentially over time. 
	
	Non-convex constrained optimization problems routinely arise in a number of areas, including robotics, machine learning, and signal processing. For instance, trajectory generation, target tracking, navigation, and path planning problems can often be cast as non-convex optimization problems \cite{schulman2014motion,malyuta2022convex, howell2019altro,bounou2023leveraging,jallet2022constrained}  with non-convex obstacle-avoidance or boundary constraints. Similarly in supervised machine learning, constraints are often imposed ensure the resulting models are both robust and fair\cite{thomdapu2023stochastic,chen2019learning}. Constraints may also be imposed to incorporate prior knowledge, e.g., to ensure that the learned features are sparse in some domain \cite{boob2020feasible}. 
	
	Existing approaches to solve \eqref{pc} include the primal-dual method \cite{jin2022stochastic}, penalty method \cite{wang2017penalty}, proximal point methods \cite{boob2024level, boob2023stochastic}, augmented Lagrangian method \cite{shi2022momentum}, and several Successive Convex Approximation (SCA) variants \cite{liu2019stochastic,liu2018online,ye2019stochastic,liu2019two,liu2021two}. The performance of these algorithms is usually measured in terms of their oracle complexity, which is the number of calls to the stochastic first-order (SFO) oracle required to reach a point that satisfies the Karush-Kuhn Tucker (KKT) conditions to within a tolerance of $\epsilon$. The primal-dual method proposed in \cite{jin2022stochastic} iteratively updates both primal and dual variables based on a stochastic approximation to augmented Lagrangian function, achieving an oracle complexity of $\mathcal{O}(\epsilon^{-2.5})$. A rate  of  $\mathcal{O}(\epsilon^{-3.5})$ is achieved by \cite{wang2017penalty} which utilizes a penalty functional to reformulate \eqref{pc} as an unconstrained optimization problem. An improved rate of $\mathcal{O}(\epsilon^{-2})$ is achieved by the  Inexact Constrained Proximal Point Algorithm (ICPPA) \cite{boob2023stochastic} and level-constrained stochastic proximal gradient (LCSPG) methods \cite{boob2024level}, both of which solve a sequence of strongly convex subproblems, obtained by adding a quadratic term to the objective and constraint functions. Both \cite{boob2023stochastic} and \cite{boob2024level} achieve $\mathcal{O}(\epsilon^{-2})$ SFO complexity using $\epsilon$-dependent step-sizes under a standard Lipschitz smoothness assumption on $U$. A related level-constrained proximal point (LCPP) method \cite{boob2020feasible} also achieved the a rate of $\O (\epsilon^{-2})$ but was customized to only handle sparsity inducing constraint functions. Improving on these works, the augmented Lagrangian method (MLALM) \cite{shi2022momentum} requires a stronger mean-squared smoothness (MSS) assumption and achieves the optimal rate of $\mathcal{O}(\epsilon^{-3/2})$, though for a fixed $\epsilon$-dependent step-size. In contrast, modern unconstrained stochastic optimization methods such as STORM and its variants have moved to using adaptive step-sizes that are generally easier to tune \cite{cutkosky2019momentum}, albeit under assumptions even stronger than MSS. 
	
	
	Different from these approaches, SCA represents a more flexible framework for solving generic non-convex constrained optimization problems. It is an iterative approach, where at every iteration, the non-convex objective and constraint functions are replaced with their convex surrogates, tailored to the problem at hand \cite{scutari2016parallel, scutari2016paralel}. Stochastic SCA variants have been proposed in \cite{yang2016parallel, liu2019stochastic, liu2018online, ye2019stochastic, liu2019two, liu2021two} for solving \eqref{pc} as well as its more general version with stochastic constraints. These works however only report asymptotic convergence results, and do not provide a comprehensive oracle complexity analysis. Non-asymptotic rates are only known for the classical SCA algorithm when solving \eqref{pc} with only convex constraints \cite{mokhtari2017large,koppel2018parallel,mokhtari2020high,idrees2021practical}, \cite{wang2024single}, where the works \cite{idrees2021practical} \cite{wang2024single} achieved the state-of-the-art rate of $\mathcal{O}(\epsilon^{-2})$. 
	
	We remark that while the oracle complexity of $\O(\epsilon^{-2})$ achieved by various algorithms that solve \eqref{pc} may be optimum under the Lipschitz smoothness assumption on $U$, it can be improved to $\mathcal{O}(\epsilon^{-3/2})$ under the MSS assumption on $f$ \cite[(4)]{arjevani2023lower}. Indeed, using momentum or variance reduction techniques \cite{levy2021storm+, tran2019hybrid, wang2019spiderboost, xu2023momentum, pham2020proxsarah}, it is well-known that the optimal oracle complexity of $\mathcal{O}(\epsilon^{-3/2})$ can be achieved for unconstrained or convex-constrained version of \eqref{pc}. Interestingly, variance reduction is applied to the level-constrained proximal algorithm in \cite{boob2024level} for finite-sum problems but the rate achieved is only $\O(\epsilon^{-2})$, with a remark that it can be improved to $\mathcal{O}(\epsilon^{-3/2})$ if $f(\x,\xib)$  is Lipschitz smooth in $\x$ for all $\xib$, 
	but without any proof \cite[Remark 5]{boob2024level}. As noted earlier, only the recently published work \cite{shi2022momentum} achieves the optimal rate of $\mathcal{O}(\epsilon^{-3/2})$ with MSS assumption, although under a restrictive condition requiring a fixed $\epsilon$-dependent step size.
	
	Before concluding this section, we remark that analysis of optimization problems with non-convex constraints inevitably requires some assumptions called the constraint qualifications (CQ). These CQs are essential regularity conditions that ensure that the constraint functions are well-behaved, so that the intermediate updates or subproblems are well-defined. Technically, these CQs are utilized to ensure that the intermediate dual variables remain bounded across iterations. Commonly used CQs include the Non-singularity condition (NSC) \cite{jin2022stochastic}, Linear Independence CQ (LICQ) \cite{lin2022complexity}, Linear Independence regularity condition (LIRC) \cite{liu2021two}, and most commonly the Mangasarian-Fromovitz CQ (MFCQ) condition \cite{ boob2023stochastic, boob2024level, boob2020feasible}. It was observed in \cite[Sec. 3.1]{boob2024level} that the classical MFCQ only implies the existence of a bound $B$ on the dual variables, but does not relate $B$ to the other problem parameters, leaving the possibility of $B$ getting arbitrarily large. It was suggested in \cite{boob2024level} to use a stronger CQ such as the strong feasibility CQ (SFCQ) so as to ensure a problem parameter-dependent oracle complexity bound. An extended version of MFCQ is used in \cite{shi2022momentum} in order to establish that the Lagrange multiplier stays bounded for all iterations. However, the extended MFCQ does not yield an explicit bound on the Lagrange multipliers, and hence the final SFO complexity in \cite{shi2022momentum} depends on an unknown (but finite) quantity $\tilde{\Lambda}$. In contrast, the slightly stronger parameterized MFCQ introduced here allows us to provide explicit bounds on the dual variables in terms of problem-dependent quantities using simpler proofs.  
	
	Table \ref{lits} summarizes the comparative performance of various state-of-the-art algorithms that can be used to solve \eqref{pc}. We observe that the near-optimal rate is achieved by the proposed CoSTA algorithm. Additional remarks on CQ, smoothness, and algorithm structure also included. 
	\subsection{Contributions} \label{sec:contri}
	The main contributions of this work are summarized below.
	\begin{itemize}
		\item \emph{Accelerated Stochastic SCA algorithm:} we incorporate the STOchastic Recursive Momentum (STORM) updates from \cite{cutkosky2019momentum, levy2021storm+} into the SCA framework to achieve near-optimal SFO complexity with non-adaptive step sizes under the MSS assumption and with adaptive step sizes under a stronger smoothness assumption.  Drawing inspiration from these works, the proposed Constrained STORM SCA (CoSTA) algorithm leverages momentum-enhanced updates to improve convergence rates.		
		\item \emph{Iteration Complexity: } non-asymptotic iteration complexity bounds are obtained by leveraging a new parameterized version of the standard MFCQ assumption, which allows us to demonstrate that the dual variables remain bounded. Together with the momentum updates, this leads to an improved near-optimal SFO complexity of $\Ot(\epsilon^{-3/2})$.
		\item \emph{Specific CoSTA algorithm design for important applications:} we demonstrate that the empirical performance of the proposed algorithm surpasses that of CSSCA \cite{liu2019stochastic}, which does not utilize momentum or variance-reduction techniques, and performs on par with LCPP \cite{boob2020feasible}, which is limited to handling sparsity-inducing constraints.
	\end{itemize}

	\renewcommand{\arraystretch}{1.35}
	\begin{table*}[h!]
		\centering
		\begin{threeparttable}[b]
			
			\caption{\label{lits}Comparison of the complexity results of several algorithms in literature to our algorithm to produce a stochastic $\epsilon$-stationary solution of constrained stochastic non-convex optimization problem (To make the comparisons fair we have converted the oracle complexity results of the all the works to match our definition of \eqref{statcond}-\eqref{complicnond2}). $L$-smoothness of $f$ implies MSS of $f$ as well as the $L$-smoothness of $U$.  $\Ot$ hides the dependence on logarithmic terms.}
			\setlength\tabcolsep{4.5pt}	
			\begin{tabular}{c|c|c|c|c|c}
				\hline
				\multicolumn{1}{c|}{\textbf{Reference}} & \multicolumn{1}{c|}{\textbf{Complexity}}                        &  
				\multicolumn{1}{c|}{{\begin{tabular}[c]{@{}c@{}}
							\textbf{Constraint} \\ \textbf{ Qualification }
				\end{tabular}}} &
				\multicolumn{1}{c|}{\textbf{Smoothness }}   &	\multicolumn{1}{c|}{{\begin{tabular}[c]{@{}c@{}}
							\textbf{Adaptive} \\ \textbf{ stepsize }
				\end{tabular}}} & \textbf{Remarks}                                                \\ \hline
				{\begin{tabular}[c]{@{}c@{}}
						CoSTA \\   \textbf{ (This work)}
				\end{tabular}}
				&   $\Ot(\epsilon^{-3/2}) $                                                                                   & Parameterized MFCQ & {\begin{tabular}[c]{@{}c@{}}
						$f$ is Smooth \\ \hline $f$ is MSS 
				\end{tabular}} &  {\begin{tabular}[c]{@{}c@{}}
						\xmark \\ \hline \cmark 
				\end{tabular}}&{\begin{tabular}[c]{@{}c@{}}
						SCA-based, non-adaptive step-size  \\   requires only MSS assumption
				\end{tabular}}               \\ \hline {\begin{tabular}[c]{@{}c@{}}
						CSSCA \cite{liu2019stochastic}, \\ SSCA \cite{ye2019stochastic}
				\end{tabular}}
				&   -                                                                                        & {\begin{tabular}[c]{@{}c@{}}
						Slater's condition established \\ for limit points
				\end{tabular}} & $f$ is Smooth &  \xmark  & Asymptotic convergence proved
				\\ \hline
				ICPPA \cite{boob2023stochastic}
				&   $\mathcal{O}(\epsilon^{-2})$                                                                                        &{\begin{tabular}[c]{@{}c@{}}
						MFCQ / SFCQ  
				\end{tabular}} & $U$ is Smooth &  \xmark   & SCA-based \\  \hline{\begin{tabular}[c]{@{}c@{}}
						Stochastic  Primal- \\ Dual SGM \cite{jin2022stochastic}
				\end{tabular}}   
				&   $\mathcal{O}(\epsilon^{-2.5})$                                                                                     &  NSC & $U$ is Smooth    &    \xmark  &  Dual variables also updated        \\ \hline  LCSPG \cite{boob2024level}
				&   $\mathcal{O}(\epsilon^{-2})$                                                                                         &{\begin{tabular}[c]{@{}c@{}}
						MFCQ / SFCQ \end{tabular}} & $U$ is Smooth  &   \xmark    & {\begin{tabular}[c]{@{}c@{}}
						Remark 5 mentions that \\
						improved rates are possible.
				\end{tabular}}            \\ \hline 
				\multicolumn{1}{c|}{LCPP\cite{boob2020feasible}}      & \multicolumn{1}{c|}{$\mathcal{O}(\epsilon^{-2})$}&  MFCQ& $U$ is Smooth  &  \xmark   & For sparsity-constrained problems only        \\ 
				\hline 
				\multicolumn{1}{c|}{MLALM\cite{shi2022momentum}}    & \multicolumn{1}{c|}{$\mathcal{O}(\epsilon^{-3/2})$}&  Extended MFCQ & $f$ is MSS  & \xmark &{\begin{tabular}[c]{@{}c@{}}
						$\O (\epsilon^{-3/2})$  but depending  \\
						on an unknown finite $\tilde{\Lambda}$ 
				\end{tabular}}        \\ \hline  
			\end{tabular}
		\end{threeparttable}
	\end{table*}
	\renewcommand{\arraystretch}{1}%
	
	\subsection{Notation and Organization} \label{notat}
	The lower case letter are used for denoting scalars while bold lower case letters are used for denoting vectors. Also, bold upper case letters are used for denoting matrices. For a vector $\a$, $\text{diag}(\a)$ is used to denote diagonal matrix with elements $(a_1,a_2, \ldots a_k)$. We collect the functions $\{g_j\}_{j=1}^J$ and $\{h_i\}_{i=1}^I$ into vector-valued functions $\g:\Rn^n\rightarrow \Rn^J$ and $\h:\Rn^n\rightarrow \Rn^I$, respectively. Binary operators such as $\leq$ and $\geq$ when applied to vectors are interpreted entry-wise, i.e., $\g(\x) \leq 0$ means $g_j(\x) \leq 0$ for all $1\leq j \leq J$. We denote $\nabla \g(\x) = \mat{\nabla g_1(\x) & \ldots \nabla g_J(\x)}^\T$.
	
	The rest of the paper is organized as follows. In Sec \ref{secc2} we present an example that motivate the problem at hand. In Section \ref{pro-al} the problem statement is posed along with the proposed algorithm. In Section \ref{con-al} a non-asymptotic convergence analysis of the CoSTA algorithm is detailed. Section \ref{app-li} briefly describes the potential application of CoSTA. Finally, Section \ref{conc} the conclusion is provided.
	
	\section{Motivating Examples} \label{secc2}
\begin{figure}
		\centering
		\includegraphics[width=0.5\textwidth]{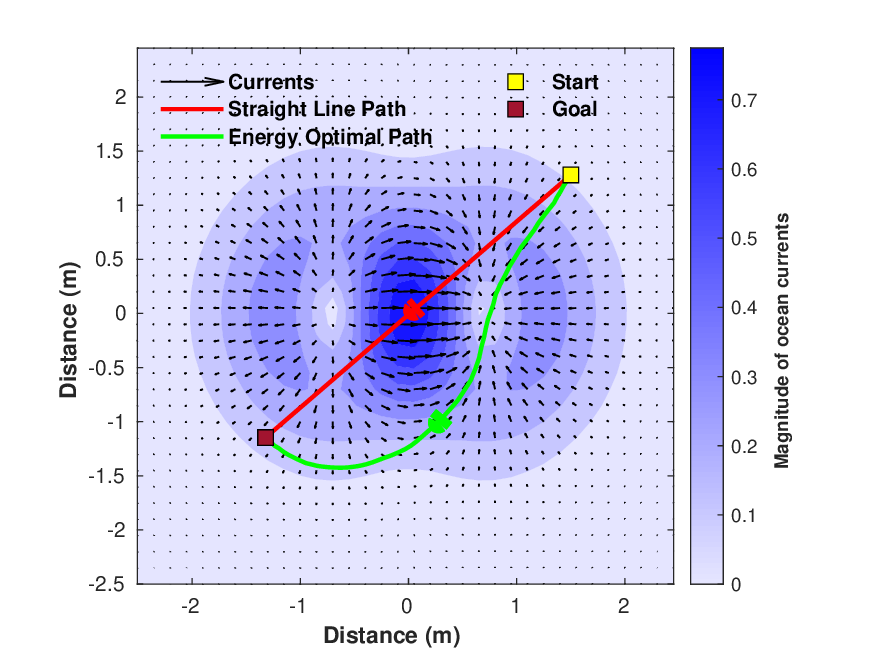}
		\caption{Navigation under ocean currents. The green energy-optimal path harnesses the currents by judiciously applying control effort as the vehicle drifts along towards the goal, achieving lower energy (8.5 J) compared to the straight line path (14.8 J).}
		\label{navig}
	\end{figure}

	In this section, we discuss a variant of the widely studied Zermelo's navigation problem \cite{Zermelo1931berDN}, and formulate it as an instance of \eqref{pc}. Consider a surface vehicle in an ocean environment that seeks to traverse from one point to another, while using minimal energy. If the ocean currents are known ahead of time, the vehicle may save energy by moving along the flows instead of moving along the straight line from the source to the destination. An example is depicted in Fig. \ref{navig}, where the straight line path (shown in red) incurs an energy of 14.8 J while the energy-optimal path (shown in green) requires only 8.5 J.  
	
	In practice however, the ocean currents in a given region $\cB\subset \Rn^2$, denoted by $\vartheta(\x)$ m/s for each $\x \in \cB$, cannot be accurately predicted. Instead, various meteorological or oceanographic agencies \cite{OceanPre5:online, HomeMerc95:online,HomeCMEM4:online,NOAATide1:online} employ ensemble forecasting, a technique that captures the measurement uncertainty inherent to the prediction process, by providing a range of possible ocean current forecasts \cite{Abetterw28:online, gneiting2005weather}. Mathematically, an instance of the ensemble forecast can be denoted by $\{\vartheta(\x,\xib)\}_{\x \in \cB}$  where $\xib$ is the index of a member of the ensemble. It has been observed that summary statistics, such as the average currents, cannot generally be used for path planning, as these may yield sub-optimal and potentially unsafe trajectories \cite{yoo2021path}. Instead, the ensemble forecasts should be directly employed to obtain trajectories that minimize the average energy. 
	
	We consider the problem of designing the energy-optimal trajectories of $N$ agents, while avoiding an obstacle. Mathematically, the problem entails finding a collection of $T$ waypoints $\{\x_i(\tau)\}_{\tau=1}^T$ for each agent $i$, collected in the super-vector $\underbar{\x} \in \Rn^{2NT}$. The waypoints are equally spaced at intervals of $\Delta t = T_f/T$ seconds where $T_f$ is the given time horizon. The $i$-th agent starts at a given location $\x_i(0) $ and wants to reach the goal location $\x_i^g$, such that $\x_i(T) = \x_i^g$. For collision avoidance, we assume that the agents and the obstacle have a circular shape with radii $r$ and $r^o$, respectively. Finally, consider the agent's movement from $\x_i(\tau)$ to $\x_i(\tau+1)$, where $0\leq \tau \leq T-1$. The control input \textemdash\ and consequently, the energy \textemdash\ can be considered proportional to the quantity $\norm{\x_i(\tau+1) - \x_i(\tau) - \vartheta(\x_i(\tau),\xib)\Delta t}^2$. Here, the term $\vartheta(\x_i(\tau),\xib)\Delta t$ approximates the vehicle's displacement without any control input. Hence, the vehicle should seek to minimize the cumulative energy input while ensuring that the magnitude of the control input does not exceed the vehicle specifications, denoted by $v_i^{\max}$. In summary, the considered problem can be written as 
	\begin{subequations}\label{eqmp}
		\begin{align}
			\min_{\underbar{\x}} &\sum_{i=1}^N  \sum_{\tau =0}^{T-1} \EE[\norm{\x_i(\tau+1) - \x_i(\tau) - \vartheta(\x_i(\tau),\xib) \Delta t}^2] \label{eqmpa} \\
			\text{s. t. } &  \x_i(T) = \x^g_i  \label{eqmp:term}\\
			&\norm{\x_i(\tau) - \x^o} \geq r^o + r, \hspace{1cm} \forall ~\tau, i \label{eqmp:nc_obs1}  \\
			&\norm{\x_i(\tau) - \x_j(\tau)} \geq 2r, \hspace{1cm} \forall ~\tau, i\neq j \label{eqmp:nc_obs2} \\
			&\hspace{-7mm} \norm{\x_i(\tau+1) - \x_i(\tau)- \vartheta(\x_i(\tau),\xib)\Delta t} \leq v_i^{\max}\Delta t \hspace{2mm} \forall ~\tau, i, \xib \label{eqmp:vr_max}
		\end{align}
	\end{subequations}
	where the objective seeks to minimize the cumulative energy consumption averaged over the given ensembles, \eqref{eqmp:term} enforces the termination condition, \eqref{eqmp:nc_obs1}-\eqref{eqmp:nc_obs2} ensure that the vehicles maintain safe distances from each other and from the obstacle, and \eqref{eqmp:vr_max} limits the maximum control input for every agent, time, and realization. Of these, \eqref{eqmp:vr_max} contains an impractically large number of constraints, one for each realization, and renders the problem difficult to solve. To avoid this issue, we consider a more stringent reformulation of the constraint, i.e., \eqref{eqmp:vr_max}. To this end, we assume that the ensemble forecasts have a uniformly bounded support, i.e., $\norm{\vartheta(\x,\xib)-\E{\vartheta(\x,\xib)}}\leq \Delta\vartheta^{\max}$ for $\x\in \cB$. Introducing $\E{\vartheta(\x_i(\tau),\xib)}\Delta t$ inside the norm on the left of \eqref{eqmp:vr_max} and  using the triangle inequality, we observe that
	\begin{align}
		&\norm{\x_i(\tau+1) - \x_i(\tau)- \vartheta(\x_i(\tau),\xib)\Delta t} \\
		& \leq \norm{\x_i(\tau+1) - \x_i(\tau) - \E{\vartheta(\x_i(\tau),\xib)}\Delta t} \nonumber\\
		&+\Delta t\norm{\vartheta(\x_i(\tau),\xib)-\E{\vartheta(\x_i(\tau),\xib)}} \nonumber \\
		&\leq \norm{\x_i(\tau+1) - \x_i(\tau) - \E{\vartheta(\x_i(\tau),\xib)}\Delta t} +\Delta  \vartheta^{\max}\Delta t. \nonumber
	\end{align}
	Hence, the constraint 
	\begin{align}
		\norm{\x_i(\tau+1) - \x_i(\tau) - \E{\vartheta(\x_i(\tau),\xib)}\Delta t} & \nonumber\\ &\hspace{-1.4cm}\leq 
		(v_i^{\max}-\Delta\vartheta^{\max})\Delta t \label{eqmp:vr_max_s}
	\end{align}
	is the tightened version of \eqref{eqmp:vr_max}. It can be seen that if we design a trajectory that satisfies \eqref{eqmp:vr_max_s}, then it would automatically satisfy \eqref{eqmp:vr_max}. However, the tightening might result in slightly sub-optimal trajectories. We can observe that the tightened problem \eqref{eqmpa}-\eqref{eqmp:nc_obs2},\eqref{eqmp:vr_max_s} admits the structure of \eqref{pc} and can be solved using the proposed algorithm.

	\section{Proposed Algorithm} \label{pro-al}
	This section describes the proposed CoSTA algorithm and states the assumptions necessary to analyze its complexity. 
	
	\subsection{Problem}
	Let us revisit \eqref{pc} to re-write it in a more compact form and introduce the associated notation. Define
	\begin{align}
		\cK \coloneqq \{\x \in \Rn^n \mid \h(\x) \leq 0\}
	\end{align}
	so that \eqref{pc} can be written as
	\begin{align}\label{pc1}
		\min_{\x \in \cK} &~~\EE[f(\x,\xib)] + u(\x) \tag{$\mathcal{P'}$}&\text{s. t. } & ~~~~ \g(\x) \leq 0.
	\end{align}
	The feasible region is given by $\cX := \cK \cap \{\x \in \Rn^n \mid \g(\x) \leq 0\}$ and the domain of \eqref{pc1} is given by 
	\begin{align}
		\cD = \text{dom}(F) \bigcap \cap_{j=1}^J \text{dom}(g_j) \bigcap \cap_{i=1}^I \text{dom}(h_i).
	\end{align}
	
	\subsection{Proposed Algorithm}
	We now describe the proposed SCA algorithm in detail. The algorithm is initialized at an arbitrary feasible point $\x_1 = \x_0 \in \cK$ that satisfies $\g(\x_1) \leq 0$. At iteration $t$ and given iterate $\x_t$, we construct a strongly convex and smooth surrogate $\fh(\x,\x_t,\xib_t)$ of the objective function $f(\x,\xib_t)$ as well as smooth and convex surrogates $\{\gt_j(\x,\x_t)\}_{j=1}^J$ (collected into the vector-valued function $\gbt(\x,\x_t)$) of the constraint functions $\{g_j(\x)\}_{j=1}^J$. These surrogates are subsequently used to formulate the following convex sub-problem that must be solved at every iteration:
	\begin{align}\label{xhata}
		\xhat =  \arg\min_{\x\in\cK}~~~ &\ft(\x,\x_t,\z_t,\xib_t)  + u(\x) \tag{$\mathcal{P}_t$} \\
		\text{s.t. } ~~~ &\gbt(\x,\x_t) \leq 0 	\nonumber
		\shortintertext{where $\ft$ is a \emph{running} approximation of $f$, given by }
		&\hspace{-3cm}\ft(\x,\x_t,\z_t,\xib_t):= \fh(\x,\x_t,\xib_t) \label{ftil}\\
		&\hspace{-2.5cm}+ (1-\beta_t)\ip{\x-\x_t}{\z_t - \nabla f(\x_{t-1},\xib_t)} \nonumber
	\end{align}
	for some $\beta_t < 1$. Here, the auxiliary variable $\z_t$ seeks to approximate $\nabla U(\x_t)$, and hence utilizes the gradient tracking update
	from \cite{cutkosky2019momentum}:
	\begin{align}
		\z_{t+1} &= \nabla f(\x_{t},\xib_{t}) + (1-\beta_{t})(\z_t - \nabla f(\x_{t-1},\xib_t)). \label{stormtrack}
	\end{align}
	We can view $\z_{t+1}$ as a convex combination of the unbiased gradient $\nabla f(\x_t,\xib_t)$ and the variance-reduced gradient estimate $\z_t + \nabla f(\x_t,\xib_t)-\nabla f(\x_{t-1},\xib_t)$ of SARAH \cite{nguyen2017sarah}, similar to the STORM update in \cite{cutkosky2019momentum} but different from the classical tracking rules used in \cite{idrees2021practical,liu2019stochastic,mokhtari2020stochastic}. Substituting \eqref{stormtrack} in \eqref{ftil}, we see that $\ft(\x,\x_t,\z_t,\xib_t)= \fh(\x,\x_t,\xib_t) + \ip{\x-\x_t}{\z_{t+1} - \nabla f(\x_t,\xib_t)}$. That is, the running approximation $\ft$ comprises of the surrogate $\fh$ and an additional linear correction term that depends on the gradient tracking error. If the surrogate follows the tangent matching property, i.e., $\nabla \fh(\x,\x,\xib)=\nabla f(\x,\xib)$, then the update in \eqref{stormtrack} also implies that 
	\begin{align}
		&\nabla \ft(\x_t,\x_t,\z_t,\xib_t) \nonumber\\
		&= \nabla \fh(\x_t,\x_t,\xib_t) +  (1-\beta_t)(\z_t - \nabla f(\x_{t-1},\xib_t)) \\
		&= \nabla f(\x_t,\xib_t) + (1-\beta_t)(\z_t - \nabla f(\x_{t-1},\xib_t)) = \z_{t+1}. \label{ze} 
	\end{align}
	Finally, we take a convex combination of $\xhat$ obtained from \eqref{xhata} and the current iterate $\x_t$ to yield the next iterate:
	\begin{align}
		\x_{t+1} &= (1-\eta_t)\x_t + \eta_t \xhat. \label{stormupdate} 
	\end{align}
	for $\eta_t < 1$. The choice of the step-size and momentum rules are similar to those in \cite{cutkosky2019momentum}. Specifically, define $G_t := \norm{\nabla f(\x_t,\xib_t)}$ so that
	\begin{align}
		\eta_t &\coloneqq \frac{\kb}{(w +\sum_{i=1}^t G_i^2)^{1/3}}, & \beta_{t+1} &= c\eta_t^2 \label{etabeta}
	\end{align}
	with $\eta_0 = \frac{\bar{k}}{w^{1/3}} < 1$. If we ensure that $\beta_t <1$, this choice adapts to the gradient norm \cite{cutkosky2019momentum} and is known to reduce the bias as well as the variance of the resulting gradient estimate \cite{tran2019hybrid}. The CoSTA algorithm is summarized in Algorithm \ref{algo1}. 
	
	\begin{algorithm}[h]
		\caption{ \textbf{Co}nstrained \textbf{ST}ORM Successive Convex \textbf{A}pproximation 
			(CoSTA) Algorithm
		} 
		\label{algo1}
		\begin{algorithmic}[1]
			\State  \textbf{Input:} Parameters $\kb,w,c$, and feasible $\x_0 = \x_1 \in \cX$, 
			\State Initialize $\beta_1  := \frac{c\kb^2}{w^{2/3}}$ and $\z_1 = \nabla f(\x_0,\xib_1) = 0$
			\For{ $t=1$ \textbf{to} $T$}  
			\State Update $\z_{t+1}$ as per \eqref{stormtrack}
			\State Evaluate $\{\eta_t, \beta_{t+1}\}$ as per \eqref{etabeta}
			\State Solve  \eqref{xhata} and update $\x_{t+1}$ as per \eqref{stormupdate}
			\EndFor
		\end{algorithmic}
	\end{algorithm}
	
	\subsection{Assumptions}\label{strsca}    
	We discuss two classes of assumptions: those regarding the problem \eqref{pc} and those regarding the choice of the surrogate functions. The first set of assumptions can be viewed as restrictions on the classes of problems that are considered here, while the second set of assumptions restrict the design choices available to us when constructing surrogate functions and applying CoSTA to a given problem. 
	
	\subsubsection{Assumptions on \eqref{pc}}
	Application of CoSTA requires \eqref{pc} to satisfy five key properties: problem should be well-defined on $\cX$, initial $\x_1$ should be feasible, $f,g_j$s and $h_i$s should be smooth and Lipschitz-continuous, variance of the objective gradient should be bounded, and finally MFCQ. 
	
	\begin{assumption}\label{cset}
		The feasible set is subsumed by the domain set $\cD$, i.e., $\cX \subseteq \cD$. 
	\end{assumption}
	
	\begin{assumption}\label{init}
		The algorithm is initialized with a feasible $\x_1$ such that $F(\x_1)-F^\star \leq B_1$.
	\end{assumption} 
	
	\begin{assumption} \label{assmooth}
		Functions $f(\cdot,\xib)$ is $L$-smooth and $G$-Lipschitz, i.e., 
		\begin{align}
			\norm{\nabla f(\x,\xib) - \nabla f(\y,\xib)} &\leq L \norm{\x-\y}, \\	
			\abs{f(\x,\xib) - f(\y,\xib)} &\leq G \norm{\x-\y} \label{fcondi}
		\end{align}
		for all $\xib$. 
	\end{assumption}
	
	\begin{assumption} \label{assmooth1}
		Functions $\{h_i\}_{i=1}^I$ are $L$-smooth and $G$-Lipschitz, i.e.,
		\begin{align}
			\norm{\nabla h_i(\x) - \nabla h_i(\y)} &\leq L \norm{\x-\y}, \\
			\abs{h_i(\x) - h_i(\y)} &\leq G \norm{\x-\y}
		\end{align}
		and functions $\{g_j\}_{j=1}^J$ are $L$-smooth and $G$-Lipschitz, 
		\begin{align}
			\norm{\nabla g_j(\x) - \nabla g_j(\y)} &\leq L \norm{\x-\y}, \\
			\abs{g_j(\x) - g_j(\y)} &\leq G \norm{\x-\y}.
		\end{align}
	\end{assumption}
	\begin{assumption}\label{asvar}
		The variances of the stochastic gradient components are bounded, i.e.,
		$\EE[\norm{\nabla f(\x,\xib) - \nabla U(\x)}^2] \leq \sigma^2$ for all $\x\in\Xc$.
	\end{assumption}	
	
	\begin{assumption}\label{mss}
		Mean Squared Smooth (MSS):
		\begin{align}
			\EE[\norm{\nabla f(\x,\xib) - \nabla f(\y,\xib)}^2] &\leq L^2 \norm{\x-\y}^2  \label{msseq1}
		\end{align}
		and Mean Squared Lipschitz (MSL)
		\begin{align}
			\EE[\abs{f(\x,\xib) - f(\y,\xib)}] &\leq G \norm{\x-\y}. \label{msseq2}
		\end{align}
	\end{assumption}
	
	\begin{assumption}\label{Amfcq}
		All feasible points $\tx \in \cX$ satisfy the parameterized $(\chi,\rho)$-MFCQ condition, implying that there exist  $\chi \geq 0$ and $\rho >0$ such that 
		\begin{subequations}\label{mfcq}
			\begin{align}
				\ip{\nabla g_j(\tx)}{\d_{\tx}} &\leq -\rho & \forall ~ j \in \cJ &:= \{k \mid g_k(\tx) \geq -\chi\}, \\
				\ip{\nabla h_i(\tx)}{\d_{\tx}} &\leq -\rho & \forall ~ i \in \cI &:= \{l \mid h_l(\tx) \geq -\chi\}, 
			\end{align}
		\end{subequations}
		for some $\d_{\tx} \in \Rn^n$ with $\norm{\d_{\tx}}_2 = 1$.
	\end{assumption} 
	
	Assumptions \ref{cset} and \ref{init} are standard and introduced to simplify the analysis. Implicit within \ref{init} is the requirement that $F^\star$ is bounded from below. The smoothness requirement in Assumption \ref{assmooth} and \ref{assmooth1} are again standard in the context of gradient-based algorithms. Note that Assumption \ref{assmooth} implies that $f$ is MSS and $U$ is $L$-smooth.  The $G$-Lipschitz condition in Assumption \ref{assmooth} implies that $G_t \leq G$ for all $t \geq 1$. Note that the MSS and MSL conditions in Assumption \ref{mss} are milder than their counterparts in Assumption \ref{assmooth}. Assumption \ref{asvar} is standard in the context of proximal stochastic gradient descent. A more relaxed version, where the gradient boundedness is only required at a specific point rather than at all points in $\cX$ has also been used in the literature \cite{khaled2023unified}, but its applicability to the non-convex constrained case remains an open problem. The MFCQ assumption in \ref{Amfcq} is critical and also common in the context of constrained optimization. The condition stems from the need to solve the convex optimization problem \eqref{xhata} at every iteration. The subsequent analysis, as well as the analysis of most widely used convex optimization solvers, necessitates that Slater's condition be satisfied for the convex subproblem, primarily to ensure that the dual variables stay bounded. The MFCQ assumption on \eqref{pc} is required to ensure that Slater's condition holds for \eqref{xhata}. 	
	
	In the literature, the classical MFCQ condition, which corresponds to $\chi = 0$ and $\rho \rightarrow 0$, is commonly encountered, but is known to yield asymptotic convergence results, as in \cite{scutari2016parallel}. This is because the classical MFCQ only implies that the dual variables associated with \eqref{xhata} are bounded by an unknown parameter that may potentially depend on the solution path, as was also observed in \cite[Remark 5]{boob2023stochastic}. Another approach, as also mentioned in Sec. \ref{sec:back}, is the extended MFCQ used in \cite{shi2022momentum} and corresponds to using $\chi = 0$ and $\rho > 0$. While the extended MFCQ guarantees a bound on the dual variables, it does not yield an explicit expression for the bound, which leaves the final SFO complexity depending on an unknown quantity $\tilde{\Lambda}$ in \cite{shi2022momentum}.  In contrast, the parameterized version of MFCQ in Assumption \ref{Amfcq} ensures that this bound is computable and depends solely on the initial conditions and the specific problem parameters such as $L$ and $\rho$. Intuitively, the vector $\d_{\tx}$ can be interpreted as a direction along which one could move a sufficiently small distance, starting from $\tx$, and stay feasible. We will formalize this intuition in Sec. \ref{con-al} and establish that Assumption \ref{Amfcq} implies the existence of a Slater point for \eqref{xhata}. The existence of such a Slater point immediately yields a bound on the dual variables. 
	
	Before concluding, we remark that one could make stronger assumptions than MFCQ that may be easier to check in practice. For instance, \cite{boob2023stochastic} established non-asymptotic results under strong feasibility and compactness of the feasible region. In the present case however, strong feasibility of \eqref{pc} does not imply that Slater's condition holds for \eqref{xhata}. Hence we do not pursue this approach. 	
	
	\subsubsection{Assumption on surrogate functions}
	We now describe the various assumptions that the surrogate functions must satisfy. 
	
	\begin{assumption}\label{Asu} For each $\y \in \cX$, the surrogate $\fh(\cdot, \y,\xib)$ satisfies:
		\begin{enumerate}
			\item Strong convexity: $\fh(\cdot, \y,\xib)$ is $\mu$-strongly convex;
			\item Smoothness: $\fh(\cdot, \y,\xib)$ is $L$-smooth; and
			\item Tangent match: $\nabla \fh(\x,\x,\xib) = \nabla f(\x,\xib)$ for all $\x \in \cX$.
		\end{enumerate}
	\end{assumption}
	
	\begin{assumption}\label{Asg}
		For each $\y \in \cX$, the surrogates $\{\gt_j(\cdot, \y)\}_{j=1}^J$ satisfy:
		\begin{enumerate}
			\item Convexity: $\{\gt_j(\cdot, \y)\}_{j=1}^J$ are convex in $\cX$;
			\item Smoothness: $\{\gt_j(\cdot, \y)\}_{j=1}^J$ are $L$-smooth in $\cX$;
			\item Lipschitz continuity: $\{\gt_j(\cdot,\y)\}_{j=1}^J$ are $G$-Lipschitz in $\cX$.
			\item Tangent match: $\nabla \gbt(\y,\y) = \nabla \g(\y)$; 
			\item Upper bound: $\g(\x) \leq \gbt(\x,\y)$ for all $\x \in \cX$ with equality at $\x =\y$. 
		\end{enumerate}
	\end{assumption}

	\begin{assumption}\label{Asul2}
		The surrogate $\ft(\cdot,\y, \z,\xib)$ is bounded over $\cX$ such that:
		\begin{align*}
			\ft(\cdot, \y,\z,\xib) + u(\x) - \min_{\x\in\cX(\y)} (\ft(\cdot, \y,\z,\xib) + u(\x)) \leq B_U
		\end{align*} 
		where $\cX(\y) := \cK \cap \{\x \mid \gbt(\x,\y) \leq 0\} \subseteq \cX$.
	\end{assumption}
	
	
	Assumptions \ref{Asu},\ref{Asg}, and \ref{Asul2} restrict the choice of the surrogate functions and have also been used in \cite{scutari2016parallel, liu2019stochastic}. Examples of surrogate functions satisfying Assumptions \ref{Asu}-\ref{Asg} can be found in \cite{sun2016majorization,scutari2016parallel} (see Sec \ref{sec:surro}). An implication of Assumptions \ref{assmooth1} and \ref{Asg} is that the function $\phi_j(\x) := g_j(\x) - \gt_j(\x,\x_t)$ is $2L$-smooth so that
	\begin{align}
		\phi_j(\x) &\leq \phi_j(\x_t) + \ip{\nabla \phi_j(\x_t)}{\x - \x_t} + L\norm{\x-\x_t}^2 \nonumber\\
		&= L\norm{\x-\x_t}^2 \label{phi1}\\
		\norm{\nabla \phi_j(\x)} &= \norm{\nabla \phi_j(\x) - \nabla \phi_j(\x_t)} \leq 2L\norm{\x -\x_t}\label{phi2}
	\end{align}
	for any $\x \in \cX$. 
	
	Assumption \eqref{Asul2} is generally satisfied if $\cX$ is compact. In that case, since $\cX(\y)\subseteq \cX$, it follows that $\cX(\y)$ is compact for all $\y$. Since $\ft(\cdot,\y,\z,\xib)$ is a smooth and convex function, it is therefore bounded over the compact set $\cX(\y)$. As an example, consider the case when $u = 0$ and let $\|\x\| \leq B_{\cX} < \infty$ for all $\x \in \cX$. Since $\fh(\cdot, \x_t, \xib_t)$ is $L$-smooth, so is $\ft(\x, \x_t, \z_t, \xib_t)$. If $\x'$ is such that $\nabla \ft(\x',\x_t,\z_t,\xib_t) = 0$, it follows from the quadratic upper bound for $L$-smooth functions that $\ft(\x,\x_t,\z_t,\xib_t) - \ft(\x',\x_t,\z_t,\xib_t) \leq \frac{L}{2}\|\x-\x'\|^2 \leq LB_{\cX}^2$. In cases where $\cX$ is not compact, it may still be possible to satisfy Assumption \eqref{Asul2} by ensuring that $\cX(\x_t)$ is compact, e.g., when $\tilde{g}_j$ is chosen to be strongly convex.

	
	
	Before concluding the discussion regarding the assumptions, the following remark regarding the problem parameters is due. 
	\begin{rem} \label{rem1}
		Given any $\chi$, it is possible to ascertain $\rho$ by solving the optimization problem:
		\begin{align} \label{rem1eq}
			\max_{\rho, \d_{\tx}} &~~\rho \\
			\text{s. t. }& \ip{\nabla g_j(\tx)}{\d_{\tx}} \leq -\rho & \forall ~~ j \in \cJ &:= \{j \mid g_j(\tx) \geq -\chi\}\nonumber\\
			&\ip{\nabla h_i(\tx)}{\d_{\tx}} \leq -\rho & \forall ~~ i \in \cI &:= \{i \mid h_i(\tx) \geq -\chi\}\nonumber \\
			& \norm{\d_{\tx}}_2 \leq 1. \nonumber
		\end{align}
		In general however, determining global values of the different problem parameters $\sigma$, $L$, $\chi$, and $\rho$ may not be easy. In practice, the algorithm parameters are often tuned directly and without first determining the problem parameters. 
	\end{rem}
	
	\subsection{Examples of Surrogate Functions} \label{sec:surro}
	The SCA framework allows surrogate functions to be customized for the specific problem at hand. Several examples of surrogates for various types of problems can be found in \cite[Sec III-A]{scutari2016parallel}. For the sake of completeness, we also mention some of them here. 
	
	For $L$-smooth constraint function $g_j$, a common choice of surrogate that satisfies all required conditions is the quadratic upper bound $\gt_j(\x,\z) = g_j(\z) + \ip{\nabla g_j(\z)}{\x-\z} + \frac{L}{2}\norm{\x - \z}^2$. If we can write $g_j(\x) = g_j^{c}(\x) + g_j^{\bar{c}}(\x)$ where the first summand is a convex function and the second summand is an $L$-smooth non-convex one, we can use the surrogate $\gt_j(\x,\z) = g_j^{c}(\x) + g_j^{\bar{c}}(\z) + \ip{\nabla g_j^{\bar{c}}(\z)}{\x-\z} + \frac{L}{2}\norm{\x - \z}^2$.
	
	In addition to these two classes of surrogates, the function $f$ admits a few more structures. For instance, if $f(\x,\xib) = f_1(f_2(\x,\xib))$ where $f_1$ is a convex function admits the surrogate $\fh(\x,\y,\xib) = f_1(f_2(\y,\xib) + \ip{\nabla f_2(\y,\xib)}{\x- \y}) + \frac{\mu}{2}\norm{\x-\y}^2$. In the case of block convex functions, the variable $\x$ can be partitioned into disjoint blocks denoted by $\{\x_k\}_{k=1}^K$ and the objective function can be expressed as $f(\x_1,\ldots,\x_k, \xib)$ where $f$ is convex in its first argument. For such functions, we can have the surrogate
	\begin{align}
			&\fh(\x,\z,\xib) = \sum_{k =1}^K \left(f(\x_k,\z_{-k},\xib) + \frac{\mu}{2}\norm{\x_k-\z_k}^2\right),  \label{scutarisurr} 
	\end{align}
	where $\z_{-k}$ collects all variables other than those in the block $\x_k$.
	\subsection{Approximate Optimality}
	The performance of the proposed algorithm will be studied in terms of its SFO complexity, i.e., the number of calls to the SFO oracle required to achieve an $\epsilon$-KKT point in expectation. Specifically, we characterize the number of SFO calls required by CoSTA to obtain a random feasible point $\tx \in \cX$, such that there exist $\lamt \in \Rn^J_+$ and $\nubt \in \Rn^I_+$ with
	\begin{align}
		&\EE\bigg[\|\nabla U(\tx) + \nabla \g(\tx)\lamt + \nabla \h(\tx)\nubt + \v_{\tx}\|\bigg] \leq \sqrt{\epsilon},  \label{statcond} \\
		& \EE\left[\lamt^\T \g(\tx)\right] \geq -\epsilon   \label{complicnond1}\\
		&\EE\left[\nubt^\T \h(\tx)\right] \geq -\epsilon\label{complicnond2}
	\end{align}
	for some $\v_{\tx} \in \partial u(\tx)$. The definition of $\epsilon$-KKT point for the stochastic case is an extension of its deterministic counterpart, first proposed in \cite{dutta2013approximate}. 
	
	\section{Complexity analysis} \label{con-al}
	In this section, we analyze the SFO complexity of the proposed algorithm. For the sake of brevity, we define for all $t \geq 1$: 
	\begin{align}
		&\ep_t = \z_{t+1} - \nabla U(\x_t)  & \text{(gradient tracking error)}, \label{et}  \\
		& \del_t  = \xhat - \x_t &\text{(iterate progress)}, \label{det} \\
		& \Delta_T = \tfrac{1}{T}\sum\nolimits_{t=1}^T\E{\norm{\del_t}} & \text{(average progress).}
	\end{align} 
	The analysis proceeds as follows. Lemmas \ref{stproof} and  \ref{le4} establish recursive bounds the $\EE[\norm{\ep_{t+1}}^2]$ and $\E{F(\x_{t+1})}$ for each $t$. Lemma \ref{feas} proves the feasibility of iterates produced by the algorithm. Theorem \ref{th1} then uses these results to establish a bound on $\Delta_T$.  Next, Lemmas \ref{slater} and \ref{dualvari} respectively establish bound on primal and dual variables required for maintaining constraint feasibility under Assumption \ref{Amfcq}. Finally, Theorem \ref{th2} uses all these results to obtain the SFO complexity of CoSTA.
	
	The results in Lemmas \ref{stproof} and \ref{le4} follow from applying the results in \cite{cutkosky2019momentum} to SCA, but using a different metric $\EE[\norm{\del_t}^2]$ instead of $\EE[\norm{\nabla U(\x_t)}^2]$. The complexity result in Theorem \ref{th1} also follows along similar lines. However, Lemmas \ref{slater} and \ref{dualvari} as well as Theorem \ref{th2} are entirely novel. The results also differ significantly from \cite{liu2019stochastic} where the focus is on non-asymptotic analysis. 
	
	\begin{lemma}\label{stproof}
		For all $t\geq 1$ and $\beta_{t+1} < \frac{1}{4}$, it holds that	
		\begin{align}
			\frac{1}{16L^2}&\E{	 \frac{\norm{\e_{t+1}}^2}{\eta_t}}  \leq \frac{1}{16L^2} \E{\frac{(1-\beta_{t+1})\norm{\e_{t}}^2}{\eta_t}}+  \frac{\E{\eta_t\norm{ \del_t}^2}}{8}  +\frac{1}{8L^2}\E{ \frac{\beta_{t+1}^2 (G_{t+1})^2}{ \eta_t }}.  \label{stbound}
		\end{align}
	\end{lemma}
	Lemma \ref{stproof} provides a bound on gradient tracking error at the $(t+1)$-th iteration in terms of the tracking error at $t$-th iteration as well as the iterate progress norm. The proof follows from expanding the left-hand side of \eqref{stbound} and introducing $ (1-\beta_{t+1})\nabla U(\x_t)$. Subsequently, using the smoothness of $f$, Assumption \ref{asvar}, and the definition of $G_{t+1}$, we obtain the desired bound. The full proof of this lemma is provided in Appendix \ref{ap1}. Next we establish the descent lemma.
	\begin{lemma} \label{le4}
		The following inequality holds for $\mu \geq \frac{ L\eta_t}{2} + \frac{3}{4}$:
		\begin{align} \label{leemm}
			\E{F(\x_{t+1}) - F(\x_t)}  \leq \frac{\E{\eta_t\norm{\e_{t}}^2}}{2}   -  \frac{\E{\eta_t\norm{\del_t}^2}}{4}. 
		\end{align}
	\end{lemma}	
	Lemma \ref{le4} bounds the change in $F$ at every iteration. The proof of Lemma \ref{le4} starts with using the smoothness of $U$ and the convexity of $u$. Subsequent use of the optimality condition of \ref{xhata} and Assumption \ref{Asu} yields the desired bound. The full proof of this lemma is provided in Appendix \ref{ap2}. Next, we show that the choice of the surrogate functions $\gt_j$ ensures that the iterates remain feasible for all $t \geq 1$. 
	
	\begin{lemma} \label{feas}
		The iterates $\{\x_t\}_{t\geq1}$ produced by Algorithm \ref{algo1} are feasible for the original problem \eqref{pc}.
	\end{lemma}
	\begin{IEEEproof}
		We use induction to establish the result. The algorithm is initialized with a feasible $\x_1$. Suppose that $\x_t \in \cX$ for some $t \geq 1$. Then, it can be seen that $\x_t$ is feasible for \eqref{xhata}, i.e., $\x_t \in \cX(\x_t)$. Further, since $\xhat$ is the solution of \eqref{xhata}, it also satisfies $\xhat \in \cX(\x_t)$. Since $\x_{t+1}$ is the convex combination of $\x_t$ and $\xhat$, and the set $\cX(\x_t)$ is convex, it follows that $\x_{t+1}\in\cX(\x_t)$. Finally, from Assumption \ref{Asg}, $\g(\x_{t+1}) \leq \gbt(\x_{t+1},\x_t) \leq 0$, implying that $\x_{t+1} \in \cX$. Therefore, by inductive hypothesis, $\x_t$ is feasible for all $t \geq 1$. 
	\end{IEEEproof}
	
	An implication of Lemma \ref{feas} is that the subproblems \eqref{xhata} are feasible. Together with Assumption \ref{Asul2}, it follows that the subproblems \eqref{xhata} are also solvable. We will now use Lemmas \ref{stproof}-\ref{feas} to establish the following theorem which bounds the average error.
	\begin{theorem} \label{th1}
		Under Assumptions \ref{cset}-\ref{asvar} and \ref{Asu},	if $0 < \bar{k} \leq w^{1/3}$, $ \tfrac{w^{2/3}}{4 \bar{k}^2}\geq c \geq 8L^2 + \frac{G^2}{3\bar{k}^3}$, $w \geq \max\{\big(\tfrac{2L\bar{k}}{4\mu -3}\big)^3,2G^2\}$, and $\mu >\tfrac{3}{4}$, then the following holds
		\begin{align}
			\Delta_T \leq \sqrt{\frac{M_T(w+TG^2)^{1/3}}{T}} = \Ot \left(\frac{1}{T^{1/3}}\right)
		\end{align}	
		where $M_T := 8 \frac{B_1}{\kb} + \frac{\sigma^2w^{1/3}}{2L^2\bar{k}^2} + \frac{c^2\bar{k}^2}{L^2}\log(T+2)$, and the $\Ot$ notation hides the dependence on logarithmic terms. 
	\end{theorem}
	
	\begin{IEEEproof}
		Since $\eta_t$ is non-increasing, the bound on $\bar{k}$ implies that $\eta_t \leq 1$ for all $t \geq 1$. We observe that since $\beta_t$ is non-increasing and $c \leq \frac{w^{2/3}}{4\kb^2}$ , it follows that $\beta_t \leq \beta_1  \leq \frac{1}{4}$  for all $t \geq 1$. Also, since $\eta_t$ is non-increasing and $w \geq \big(\tfrac{2L\bar{k}}{4\mu -3}\big)^3$, it follows that $\mu \geq \frac{L\eta_0}{2}+\frac{3}{4} \geq \frac{L\eta_t}{2}+\frac{3}{4}$, which is the required condition in Lemma \ref{le4}.
		
		Defining $\Psi_t = \E{F(\x_t)-F^\star} + \tfrac{1}{16L^2}\EE[\tfrac{\norm{\e_{t}}^2}{  \eta_{t-1}}]$ and adding \eqref{stbound} with \eqref{leemm}, we obtain
		\begin{align}
			\Psi_{t+1}-\Psi_t &\leq \E{\frac{\alpha_t}{16L^2}\norm{\e_{t}}^2} -\frac{\E{\eta_t\norm{\del_t}^2}}{8} +  \frac{c^2}{8L^2}  \E{\eta_t^3 G_{t+1}^2}, \label{stproof1}
		\end{align}
		where $\alpha_t := 8L^2\eta_t + \frac{(1-\beta_{t+1})}{\eta_t} - \frac{1}{\eta_{t-1}}$. In \eqref{stproof1}, we observe that the second term is already negative. We also want to select the algorithm parameters such that the first term is also non-positive, leaving the rate depending on the third term, which we will subsequently bound. 
		
		We first show that $\alpha_t \leq 0$. We begin with observing from concavity that since $(a + b)^{1/3}-a^{1/3} \leq \tfrac{b}{3a^{2/3}}$ for all $a,b \geq 0$, we have that 
		\begin{align}
			\frac{\bar{k}}{\eta_t}-\frac{\bar{k}}{\eta_{t-1}} &= (w+ \sum_{i =1}^{t}G_i^2)^{1/3} - (w+ \sum_{i =1}^{t-1}G_i^2)^{1/3}  \nonumber\\
			&\leq \frac{G_t^2}{3(w+ \sum_{i =1}^{t-1}G_i^2)^{2/3}} = \frac{G_t^2\eta_{t-1}^2}{3\bar{k}^2}. \label{th1proof1}
		\end{align}
		Therefore, we can simplify the last two terms in the definition of $\alpha_t$ by using \eqref{th1proof1} and the definition of $\beta_t$ as
		\begin{align}
			\frac{(1-\beta_{t+1})}{\eta_t} - \frac{1}{\eta_{t-1}} &= \frac{1}{\eta_t} - \frac{1}{\eta_{t-1}} - \frac{\beta_{t+1}}{\eta_t} \label{th1proof2}\\
			&\leq \frac{G_t^2\eta_{t-1}^2}{3\bar{k}^3} - c\eta_t. \label{th1proof4}
		\end{align}
		Finally, since $\eta_t \leq \eta_{t-1}$, $\beta_t > 0$, and $G_t \leq G$, we have that
		\begin{align}
			\alpha_t &=  \frac{(1-\beta_{t+1})}{\eta_t} - \frac{1}{\eta_{t-1}} + 8L^2\eta_t \\
			&\leqtext{\eqref{th1proof4}} \frac{G^2\eta_{t-1}^2}{3\bar{k}^3} - c\eta_{t-1} + 8L^2\eta_{t-1} \\
			&= -2\eta_{t-1}\left(\frac{c}{2} - 4L^2 - \frac{G^2\eta_{t-1}}{6\bar{k}^3}\right)  \leq -2L^2d\eta_{t-1}  \label{th1proof3}
		\end{align}
		where $d:= \tfrac{1}{2L^2}(c - 8L^2 - \tfrac{G^2}{3\bar{k}^3}) > 0$. Therefore the first term on the right of \eqref{stproof1} can be dropped. 
		
		For the second term on the right of \eqref{stproof1}, using the fact that $\eta_T \leq \eta_t$ as well as the Cauchy-Schwarz inequalities ($\E{\sX\sY}^2\leq \E{\sX^2}\E{\sY^2}$ for random variables $\sX,\sY$ and $\mathbf{1}^\T\w \leq \sqrt{T}\norm{\w}$ for vector $\w \in \Rn^T$), we obtain
		\begin{align}
			&\E{\sum_{t=1}^T\eta_t\norm{\del_t}^2} \geq	\E{\sum_{t=1}^T\eta_T\norm{\del_t}^2} \label{delta1} \\
			&\geq \frac{1}{\E{1/\eta_T}}\E{\sqrt{\sum_{t=1}^T\norm{\del_t}^2}}^2 \geq \frac{\bar{k}T\Delta_T^2}{(w+TG^2)^{1/3}}, \label{delta}
		\end{align}
		where we have used the fact that $\tfrac{1}{\eta_T} \leq \tfrac{1}{\bar{k}}(w+TG^2)^{1/3}$.
		
		For the third term on the right of \eqref{stproof1}, we can use the fact that $w \geq 2G^2$ and the concavity of the log function as in \cite[Lemma 4]{cutkosky2019momentum} to obtain
		\begin{align}
			&\sum\nolimits_{t=1}^T\eta_t^3 G_{t+1}^2 =  \kb^3   \sum\nolimits_{t=1}^T \frac{ G_{t+1}^2}{(w+\sum_{i=1}^{t} G_i^2)}  \nonumber \\
			&\leq \kb^3   \sum\nolimits_{t=1}^T \frac{ G_{t+1}^2}{(G^2+\sum_{i=1}^{t+1} G_i^2)}  \leq  \kb^3   \ln \Big(1 + \sum\nolimits_{t=1}^{T+1} \frac{G_t^2}{G^2}\Big) \nonumber \\
			& \leq \kb^3 \ln (T + 2). \label{logsum}
		\end{align}
		Substituting \eqref{delta} and \eqref{logsum} into \eqref{stproof1}, using the fact that $\beta_t \leq 1/4$, and summing over $t = 1, 2, \ldots, T$, we obtain
		\begin{align}
			\Psi_{T+1}-\Psi_1 &\leq -\frac{\bar{k}T\Delta_T^2}{8(w+TG^2)^{1/3}} + \frac{c^2\bar{k}^3}{8L^2} \ln(T+2). \label{psisum}
		\end{align}
		Here, since $\EE[\norm{\e_1}^2] \leq \sigma^2$ from Assumption \ref{asvar}, we can bound
		\begin{align}
			\Psi_1 \leq F(\x_1) - F^\star + \frac{\sigma^2}{16L^2\eta_0}. \label{psibound}
		\end{align}
		Further, the feasibility of $\x_{T+1}$ as shown in Lemma \ref{feas} implies that $F^\star \leq F(\x_{T+1})$ and hence $\Psi_{T+1}\geq 0$. Substituting into \eqref{psisum} and re-arranging, we obtain
		\begin{align}
			\Delta_T^2 \leq M_T\frac{(w+TG^2)^{1/3}}{T} 
		\end{align}
		where $M_T$ is as defined in the statement of the theorem. 
	\end{IEEEproof}
	\begin{corollary} \label{cor1}
		Under Assumptions \ref{cset}-\ref{asvar} and \ref{Asu}, along with the conditions of Theorem \ref{th1}, it also holds that
		\begin{subequations}
			\begin{align}
				\tfrac{1}{T} \sum_{t =1}^{T}  \E{\norm{\del_t}^2}  \leq \frac{M_T(w+TG^2)^{1/3}}{T} = \Ot \left(\frac{1}{T^{2/3}}\right) \label{cor:eq1} \\
				\tfrac{1}{T}\sum_{t=1}^T\E{\norm{\e_t}} \leq \sqrt{\frac{M_T(w+TG^2)^{1/3}}{dT}} = \Ot \left(\frac{1}{T^{1/3}}\right) \label{cor:eq2}
			\end{align}	
		\end{subequations}
		where $M_T := 8 \frac{B_1}{\kb} + \frac{\sigma^2w^{1/3}}{2L^2\bar{k}^2} + \frac{c^2\bar{k}^2}{L^2}\log(T+2)$, 
		and $d= \tfrac{1}{2L^2}(c - 8L^2 - \tfrac{G^2}{3\bar{k}^3})$.
	\end{corollary}
	The bound in \eqref{cor:eq1} of Corollary \ref{cor1} can be established along the lines of Theorem \ref{th1} by lower bounding the RHS of \eqref{delta1} and using $\eta_T ~\geq ~	\tfrac{\kb}{(w+TG^2)^{1/3}}$. The bound in  \eqref{cor:eq2} of Corollary \ref{cor1} can be established on similar lines. Here, instead of dropping the first term on the right of \eqref{stproof1}, we drop the second term instead, and obtain a bound on the first term by proceeding along similar lines. 
	
	It is important to note that in Lemma \ref{stproof}, the last term of \eqref{stbound} affects the upper bounds, causing the convergence rate to depend on a logarithmic term. However, if the adaptive step sizes in \eqref{etabeta} are replaced with non-adaptive choices, the bounds in Theorem~\ref{th1} and Corollary~\ref{cor1} can be recovered under a milder assumption, as formalized in the following corollary. The detailed proof is provided in Appendix~\ref{prc2c3}.
	\begin{corollary} \label{cor2}
			Using non-adaptive step size, namely $\eta_t = \tfrac{\kb}{(w + t)^{1/3}}$ and $\beta_{t+1} = c\eta_t^2$, the bounds stated in Theorem~\ref{th1} and Corollary~\ref{cor1} continue to hold under the weaker MSS Assumption~\ref{mss}, with slightly modified constant coefficients.
	\end{corollary}
	Theorem \ref{th1} and Corollary \ref{cor1} establish bounds on the average tracking error and the average progress. These results will later be used to establish the required $\epsilon$-KKT conditions and hence the SFO-complexity bound. To establish these bounds however, we need a strict feasibility condition, which will help us bound the dual variables corresponding to the $t$-th sub-problem \eqref{xhata}. This is where the MFCQ condition in Assumption \ref{Amfcq} enters into the picture. 
	\begin{lemma}\label{slater}
		Under Assumption \ref{Amfcq}, there exists a strictly feasible $\tx(\x_t) \in \cX$ such that
		\begin{align}
			\gt_j(\tx(\x_t),\x_t) &\leq -\frac{\rho^2}{2L}, &j= 1, . . . , J. \label{gbt}\\
			h_i(\tx(\x_t)) &\leq -\frac{\rho^2}{2L}, &i = 1, . . . , I. \label{hval}
		\end{align}
		where $\rho$ is the parameter in Assumption \eqref{Amfcq} for $\chi \geq \tfrac{\rho}{L}(G+\rho/2)$.
	\end{lemma}
	Lemma \ref{slater} implies that a strong version of Slater's constraint qualification holds for \eqref{xhata}. In general, one cannot directly assume the existence of such as Slater point because the intermediate subproblems depend on $\x_t$ which is generated by the algorithm. The proof of Lemma \ref{slater} is provided in Appendix \ref{ap3} and follows  from the smoothness and Lipschitz continuity properties of $\gt_j$ and $h_i$.
	
	Slater's constraint qualification together with the fact that $F^\star > -\infty$ implies that strong duality holds for \eqref{xhata}. Therefore, the primal-dual optimum solution, denoted by $(\xhat, \lamh_t, \nubh_t)$, satisfies the KKT conditions:
	\begin{subequations}\label{kkt}
		\begin{align}
			\gbt(\xhat,\x_t) &\leq 0, \h(\xhat) \leq 0 \label{pfeas}\\
			\lamh_t &\geq 0, \nubh_t \geq 0  \label{dfeas}\\
			\lamh_t \odot \gbt(\xhat,\x_t) &= \nubh_t \odot \h(\xhat)= 0\label{cs}
		\end{align}
		where $\odot$ denotes the Hadamard (entry-wise) product. Further, for some $\v_{\xhat} \in \partial u(\xhat)$, it holds that
		\begin{align}
			& \nabla \ft(\xhat,\x_t,\z_t,\xib_t)+ \v_{\xhat} + \nabla\gbt(\xhat,\x_t)\lamh_t +\nabla \h(\xhat)\nubh_t  = 0. \label{stat2}
		\end{align}
	\end{subequations}
	The Slater's condition along with the boundedness assumption in \ref{Asul2} allows us to bound the dual variables for each \eqref{xhata}. 
	%
	\begin{lemma} \label{dualvari}
		Under conditions of Lemma \ref{slater}, for $t \geq 1$, we have that $\|\lamh_t\|_1 + \|\nubh_t\|_1 \leq \frac{2B_U L}{\rho^2}$.
	\end{lemma}
	
	To prove Lemma \ref{dualvari} we use the convexity of $\ft$, $\gt_j$, $h_i$, and $u$, as well as the results in Lemma \ref{slater} and Assumption \ref{Asul2}. The proof of Lemma \ref{dualvari} is provided in Appendix \ref{ap4}. Finally, we state the main result, which is regarding the existence of an $\epsilon$-optimal KKT point. 
	
	
	\begin{theorem} \label{th2}
		Under Assumptions \ref{cset}-\ref{asvar} and \ref{Amfcq}-\ref{Asul2}, along with the conditions of Theorem \ref{th1}, Corollary \ref{cor1} and Lemma \ref{slater}, there exists some $1\leq t \leq \Ot(\epsilon^{-3/2})$ and $\lamh_t \in \Rn^I_+$, $\nubh_t \in \Rn^J_+$ such that the point $(\xhat, \lamh_t, \nubh_t)$ is $\epsilon$-KKT optimal. 
	\end{theorem}
	Equivalently, Theorem \ref{th2} states that the point $(\xhat, \lamh_t, \nubh_t)$ for some $1\leq t \leq T$  
	is $\epsilon$-KKT optimal for $\epsilon = \Ot \left(T^{-2/3}\right)$. We provide the proof of Theorem \ref{th2}, that relies on the results established in Theorem \ref{th1} and Lemma \ref{dualvari}.
	
	\begin{IEEEproof}
		We with begin with establishing that $(\xhat,\lamh_t,\nubh_t)$ satisfies the approximate stationarity condition in \eqref{statcond}  for some $1\leq t\leq T$. Observe from \eqref{ftil} and \eqref{ze} that
		\begin{align}
			\nabla &\ft(\xhat,\x_t,\z_t,\xib_t) \nonumber\\
			&= \nabla \fh(\xhat,\x_t, \xib_t) + (1-\beta_t)(\z_t - \nabla f(\x_{t-1},\xib_t)) \\
			&=\nabla \fh(\xhat,\x_t, \xib_t) - \nabla \fh(\x_t,\x_t, \xib_t) + \z_{t+1}.
		\end{align}
		Substituting in \eqref{stat2} and re-arranging, we obtain
		\begin{align} \label{prf11}
			\v_{\xhat} &+ \nabla \h(\xhat)\nubh_t = - \nabla \gbt(\xhat,\x_t)\lamh_t  -\nabla\fh(\xhat,\x_t, \xib_t) + \nabla\fh(\x_t,\x_t, \xib_t) - \z_{t+1}.
		\end{align}
		Adding $\nabla U(\xhat) + \nabla \g(\xhat)\lamh_t$ to both sides, we obtain
		\begin{align}
			&\nabla U(\xhat) + \nabla \g(\xhat)\lamh_t  + \nabla \h(\xhat)\nubh_t +\v_{\xhat} \nonumber\\
			& = \nabla U(\xhat) - \z_{t+1} + \nabla\fh(\x_t,\x_t, \xib_t)-\nabla\fh(\xhat,\x_t, \xib_t)\nonumber\\
			& +  \Big ( \nabla \g(\xhat) - \nabla \gbt(\xhat,\x_t) \Big )\lamh_t =: \pib_t.
		\end{align}
		The norm of $\pib_t$ can now be bounded using the triangle inequality. First note that since $U$ is $L$-smooth, we have that
		\begin{align}
			\norm{\nabla U(\xhat) - \z_{t+1}} &= \norm{\nabla U(\xhat) - \nabla U(\x_t) + \nabla U(\x_t) - \z_{t+1}} \nonumber\\
			&\hspace{-1cm}\leq L\norm{\xhat-\x_t} + \norm{\ep_t} = L\norm{\del_t}+\norm{\ep_t}.
		\end{align}
		Likewise, since $\fh(\cdot,\x_t,\xib_t)$ is $L$-smooth from Assumption \ref{Asu}, it follows that $\norm{\nabla\fh(\x_t,\x_t, \xib_t)-\nabla\fh(\xhat,\x_t, \xib_t)}\leq L\norm{\del_t}$. Also from \eqref{phi2} and since $\lamh_t\geq 0$, we have that 
		\begin{align}
			&\norm{[\lamh_t]_j(\nabla g_j(\xhat) - \nabla \gt_j(\xhat,\x_t))}  \leq 2L [\lamh_t]_j \norm{\xhat-\x_t}
		\end{align}
		and hence from the triangle inequality,
		\begin{align}
			&\norm{\Big ( \nabla \g(\xhat) - \nabla \gbt(\xhat,\x_t) \Big )\lamh_t} \nonumber\\
			&\leq \sum_{j=1}^J \norm{[\lamh_t]_j(\nabla g_j(\xhat) - \nabla \gt_j(\xhat,\x_t))} \nonumber\\
			&\leq 2L\norm{\lamh_t}_1\norm{\xhat-\x_t}.
		\end{align}
		Combining these bounds, we obtain
		\begin{align}
			\norm{\pib_t} &\leq 2L(1+\tfrac{2B_UL}{\rho^2})\norm{\del_t} + \norm{\ep_t}.
		\end{align}
		Taking expectation and summing over all $t = 1, 2, \ldots, T$, and using the bounds in Theorem \ref{th1} and Corollary \ref{cor1}, we obtain 
		\begin{align}
			&\tfrac{1}{T}\sum\nolimits_{t=1}^T \E{\norm{\pib_t}}^2 \nonumber \\
			&~~~~ \leq \tfrac{2}{T}\sum_{t=1}^T \left(2L + \frac{4B_UL^2}{\rho^2}\right)^2  \E{\norm{\del_t}}^2+ \tfrac{2}{T}\sum_{t=1}^T \E{\norm{\ep_t}}^2    \label{frate}  \\
			&~~~~ \leq \left( M_T \left(2L + \frac{4B_UL^2}{\rho^2}\right)^2 +\frac{M_T}{d} \right)\frac{2 (w+TG^2)^{1/3}}{T} \nonumber\\
			&~~~~ =\Ot(T^{-2/3}). \label{st1}
		\end{align} 
		Also note that,
		\begin{align}
			&\lamh_t^\T\g(\xhat) = \sum\nolimits_{j=1}^J[\lamh_t]_j\big(\gt_j(\xhat,\x_t)+ g_j(\xhat) - \gt_j(\xhat,\x_t)\big) \nonumber \\
			&\eqtext{\eqref{cs}}-\sum\nolimits_{j=1}^J[\lamh_t]_j\big(\gt_j(\xhat,\x_t) - g_j(\xhat) \big)  \\
			&\geqtext{\eqref{phi1}}  -L \norm{\del_t}^2 \sum\nolimits_{j=1}^J[\lamh_t]_j \geq -\frac{2B_U L^2}{\rho^2}\norm{\del_t}^2.  \label{complll}
		\end{align}
		Taking expectation, summing over $t = 1, 2, \ldots, T$, and using the bound in Corollary \ref{cor1} we obtain,
		\begin{align}
			\tfrac{1}{T}\sum\nolimits_{t =1}^{T} \EE [\lamh_t^\T\g(\xhat)] &\geq -\frac{2B_U L^2}{\rho^2} \frac{M_T(w+TG^2)^{1/3}}{T}\nonumber \\
			& = \Ot (T^{-2/3}). \label{cs1}
		\end{align}
		Adding \eqref{st1} with \eqref{cs1} and noting that since $\min_{1\leq t\leq T} (\E{\norm{\pib_t}}^2 - \EE [\lamh_t^\T\g(\xhat)]) \leq \tfrac{1}{T}\sum\nolimits_{t=1}^T \E{\norm{\pib_t}}^2 -   \tfrac{1}{T}\sum\nolimits_{t =1}^{T} \EE [\lamh_t^\T\g(\xhat)]  $,	it follows that there exists $t\in\{1, \ldots, T\}$ such that $\E{\norm{\pib_t}} \leq \sqrt{\epsilon}$ and  $\E{\lamh_t^\T\g(\xhat)} \geq - \epsilon $.
		Comparing with the required condition in \eqref{statcond} and \eqref{complicnond1}, we can set $\epsilon = T^{-2/3}$. In other words, there exists some $t$ such that the point $(\xhat,\lamh_t,\nubh_t)$ is approximately stationary and satisfies approximate complementary slackness.
		
		For the same $t$, it already holds that $\lamh_t \geq 0, \nubh_t \geq 0$ from \eqref{dfeas}. Further, $\xhat$ is also feasible for \eqref{pc} from \eqref{pfeas}  and since $g_j(\xhat) \leq \gt_j(\xhat,\x_t) \leq 0$. Likewise, \eqref{complicnond2} already holds from  \eqref{cs}.
	\end{IEEEproof}
	
	Under a non-adaptive step-size setting with $\eta_t = \frac{\kb}{(w+t)^{1/3}}$ and $\beta_{t+1} = c\eta_t^2$, we state the following corollary, which follows along the lines of Theorem \ref{th2} and whose proof is detailed in Appendix \ref{prc2c3}.
	\begin{corollary} \label{cor3}
			Under milder MSS assumption \ref{mss} and using a non-adaptive step size, along with the bounds established in Corollary \ref{cor2} and Lemma \ref{slater}, there exists some  $1 \leq t \leq \Ot(\epsilon^{-3/2})$ and $\lamh_t \in \Rn^I_+$, $\nubh_t \in \Rn^J_+$ such that the point $(\xhat, \lamh_t, \nubh_t)$ is $\epsilon$-KKT optimal. 	
	\end{corollary}

Under the adaptive step-size assumption, the SFO complexity bound in Theorem \ref{th2} matches that of STORM \cite{cutkosky2019momentum} and nearly matches Pstorm \cite{xu2023momentum},  STORM+ \cite{levy2021storm+} achieving a complexity of $\mathcal{O} (\epsilon^{-3/2})$ for the unconstrained case. Notably, while STORM and Pstorm were proposed for unconstrained and convex-constrained problems, respectively, CoSTA extends this to non-convex constrained framework with a comparable SFO complexity under similar assumptions. Furthermore, under the non-adaptive step-size setting, CoSTA attains the same convergence rate but under milder MSS assumption. A key distinction of CoSTA, however, lies in its ability to achieve a near-optimal rate (i.e., $\Ot (\epsilon^{-3/2})$) under adaptive step sizes, a feature that distinguishes it from \cite{shi2022momentum}, which employs non-adaptive step sizes. It is remarked that the non-asymptotic results in Theorem \ref{th2} and Corollary \ref{cor3} can also be used to infer the asymptotic convergence of the corresponding quantities.

\section{Simulations} \label{app-li}
In this section we will evaluate the numerical performance of CoSTA on the two problems. 
We will first discuss the classical problem of sparse binary classification and compare the performance of CoSTA with the Level-Constrained Proximal Point (LCPP) algorithm designed to handle non-convex sparse models \cite{boob2020feasible}. Subsequently, we will discuss the implementation of CoSTA on the energy-optimal trajectory planning problem discussed in Sec. \ref{secc2}.


\subsection{Sparse Binary Classification} \label{sec:SBC}
The problem of finding a sparse binary classifier can be written in the form of \eqref{pc1} with the objective being the logistic loss function, i.e., $f(\x,\{\a,b\}) = \log (1+ e^{- b \a^\T\x})$ and the constraint is given by $\bar{g}(\x) \leq \tau$, where $\bar{g}$ is a sparsity-inducing function such as the $\ell_1$ norm or the minimax concave penalty (MCP) \cite{boob2020feasible}. Here, the random variable denotes independent observations $\{\a, b\}$. In \cite{boob2020feasible}, the constraint function $\bar{g}$ is chosen to be the MCP, defined as: $\bar{g}(\x) = \lambda \|\x\|_1 - \sum\nolimits_{k=1}^{n} h_{\lambda,\theta}([\x]_k) $
where $[\x]_k$ denotes the $k$-th entry of $\x$ and  $h_{\lambda,\theta} $ is defined as:
\begin{align} \label{mcpp}
	h_{\lambda,\theta}(x) = \begin{cases}
		\frac{x^2}{2\theta}, &  |x|\leq \theta \lambda  \\
		\lambda|x|- \frac{\theta\lambda^2}{2}, &  |x| > \theta \lambda. \\
	\end{cases}
\end{align}
Since $\bar{g}$ in \eqref{mcpp} is non-smooth, we consider an alternative smooth constraint given by,  
\begin{align} \label{contt}
	g(\x) & =\sum\nolimits_{k=1}^{n} \Big(h_{\lambda,\varrho}([\x]_k)  -  h_{\lambda,\theta}([\x]_k)  \Big)  \leq \tau
\end{align}
\noindent where $\varrho, \theta, \lambda >0 $ and $\varrho \ll \theta$. This results in a smooth, non-convex function $g$ that can be explicitly written as
\begin{align} 
	&g(x) = \begin{cases}
		\frac{x^2}{2} \Big(\frac{1}{\varrho} -  \frac{1}{\theta}\Big)  ,  &  |x|\leq \varrho \lambda  \\
		\lambda|x|- \frac{\varrho\lambda^2}{2}- \frac{x^2}{2\theta}, & \varrho
		\lambda  < |x| \leq  \theta \lambda \\
		\frac{(\theta -\varrho)\lambda^2}{2} , & |x| > \theta \lambda.
	\end{cases}   
	\end{align}
	Since the objective is already convex, its running surrogate $\ft$ can be constructed using $\fh(\x,\x_t,\xib_t) = f(\x_t,\xib_t) + \ip{\nabla f(\x_t,\xib_t)}{\x -\x_t} + \frac{\mu}{2} \norm{\x - \x_t}^2$ and \eqref{ftil} as
	\begin{align}
&\ft(\x,\x_t,\z_t,\xib_t)=  f(\x_t,\xib_t) + \langle \nabla f(\x_t,\xib_t) \\
&+ (1-\beta_t)(\z_t - \nabla f(\x_{t-1},\xib_t)), \x- \x_t \rangle + \frac{\mu}{2} \norm{\x - \x_t}^2 \nonumber 
\end{align}
where $\xib_t = \{\a_t,b_t\}$. The constraint $g$ in \eqref{contt} is non-convex and hence, the surrogate $\gt$ at point $\x_t$ can be constructed as,

\begin{align}  
	&\gt(\x,\x_t)    \label{surroeg} \\ 
	&= \sum\nolimits_{k=1}^{n} \Big( \tilde{g}_{lp}([\x]_k,[\x_t]_k)+ \tfrac{L}{2}([\x]_k - [\x_t]_k )^2\Big) \leq \tau  \nonumber
\end{align}
where $\tilde{g}_{lp}$ is defined as
\begin{align}
	&	\tilde{g}_{lp}(x, x_t)  = g(x_t) + \ip{\nabla g(x_t)}{ x- x_t} \nonumber \\
	&=   
	\begin{cases}
		\frac{x_t x}{\varrho} - \frac{x_t x}{\theta}  -  \frac{x_t^2}{2\varrho}  +   \frac{x_t^2}{2\theta},  &  |x_{t}|\leq \varrho \lambda  \\
		\lambda \cdot x\text{sign}(x_t)- \frac{x_t x}{\theta} - \frac{\varrho\lambda^2}{2}+ \frac{x_t^2}{2\theta} 	, &  \varrho\lambda  < |x_{t}| \leq  \theta \lambda \\
		\frac{(\theta -\varrho)\lambda^2}{2}, & |x_{t}| > \theta \lambda.
	\end{cases}  \label{surroused}
	\end{align}
	Observe that \ref{surroeg} satisfies \ref{Asg}. We compare the performance of our algorithm with that of the LCPP method, which was proposed in \cite{boob2020feasible} for solving stochastic non-convex optimization problems with sparsity inducing non-convex constraints. It is proximal point algorithm that solves the a sequence of convex subproblems with gradually relaxed constraint levels. 
	
	\begin{table}[h]
\footnotesize
\centering
\caption{\label{dataset} Dataset Description}
\begin{tabular}{c|c|c|c}
	\hline
	\footnotesize
	Datasets & Training size & Test Size & Dimentionality ($n$)\\ \hline
	MNIST    & $60000$         & $10000$     & $784$            \\  \hline
	gisette  & $6000$          & $1000  $    & $5000$           \\ 
	\hline
\end{tabular}
\end{table}

For our simulations use two dfatasets: MNIST and gisette, details of which are summarized in Table \ref{dataset}. For the MNIST dataset, we consider the binary classification problem of classifying the handwritten digit 5 from the other digits. As the dataset consists of 10 classes this can be achieved by assigning a label $+1$ to digit $5$ and $-1$ to other digits. 
Likewise, gisette dataset also consists handwritten digits, but here the classification is between the often confused digits 4 and 9. 

For MCP we take  $\lambda = 2 $ and $\theta =  5$, which correspond to one of the parameter values used in \cite{boob2020feasible}. The following CoSTA parameters were used after tuning them for best performance: (i) $\mu =  0.05, \kb = 0.0051, c = 6 \times 10^{5}, w = 9000, \tau = 0.05n$ for gisette; and (ii) $\mu =  0.06, \kb = 0.0018, c = 1.4 \times 10^{6}, w = 38000, \tau = 0.1n$ for MNIST.

\begin{figure}
\label{perrr} 
\centering
\begin{subfigure}[b]{0.53\linewidth}
	\centering			\includegraphics[width=\linewidth, trim={3cm 0 0.7cm 0},clip]{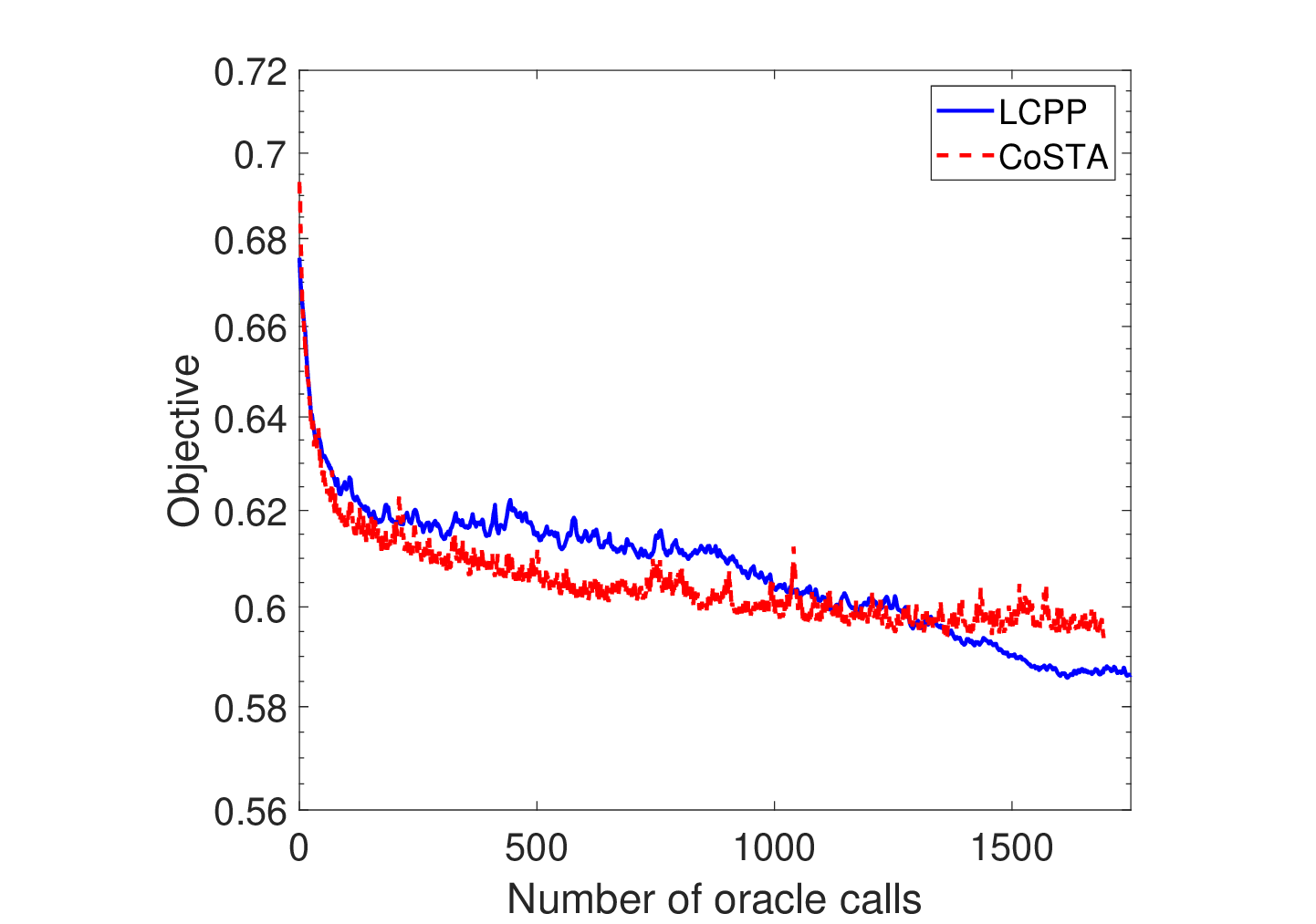}  
	\caption{MNIST}
	\label{perf}
\end{subfigure} 
\hspace{-1cm}
\hfill
\hspace{-1cm}
\begin{subfigure}[b]{0.53\linewidth}
	\centering
	\includegraphics[width=\linewidth, trim={3cm 0 0.7cm 0},clip]{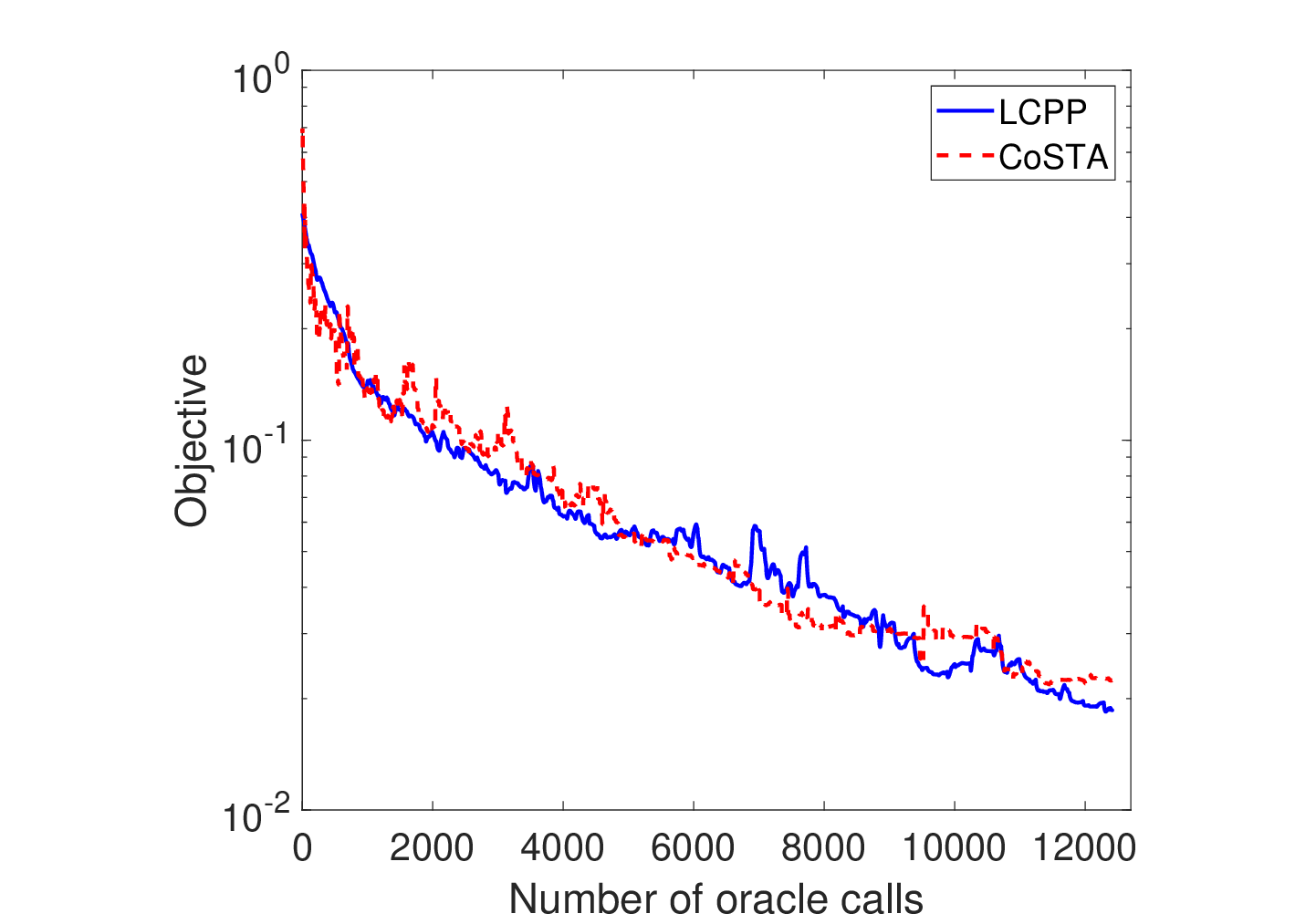}  
	\caption{gisette}	
	\label{perffm}
\end{subfigure} 
\caption{Evolution of the objective function value across the number of calls to the oracle. }
\end{figure}

Figs. \ref{perf} and \ref{perffm} show the performance of CoSTA and LCPP for the MNIST and gisette datasets, respectively. We observe that the proposed CoSTA algorithm performs at par with LCPP for both the datasets. The classification performance is also summarized in Table \ref{ress}, where it can be seen that the test performance of CoSTA is always better than that of LCPP. These results are remarkable since CoSTA is a general-purpose non-convex optimization algorithm while LCPP is a customized algorithm for working with non-convex sparse models. 

\begin{figure*}[]
\centering
\begin{subfigure}[b]{0.3468\textwidth}
	\centering
	\includegraphics[width=\textwidth, trim={0cm 0 3.2cm 0},clip]{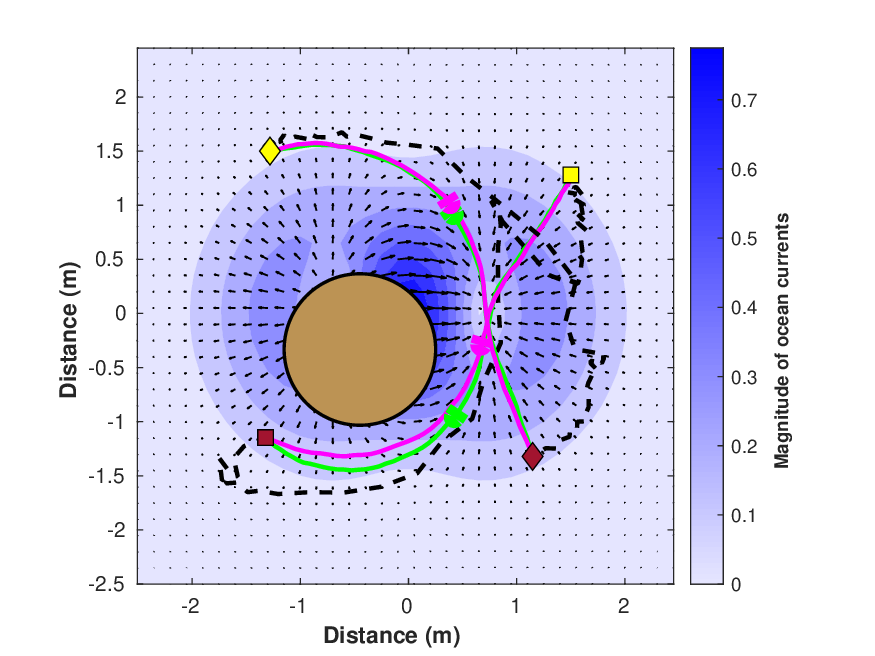}
\end{subfigure}
\hspace{-0.7cm}
\hfill
\hspace{-0.7cm}
\begin{subfigure}[b]{0.302\textwidth}
	\centering
	\includegraphics[width=\textwidth, trim={1.5cm 0 3.2cm 0},clip]{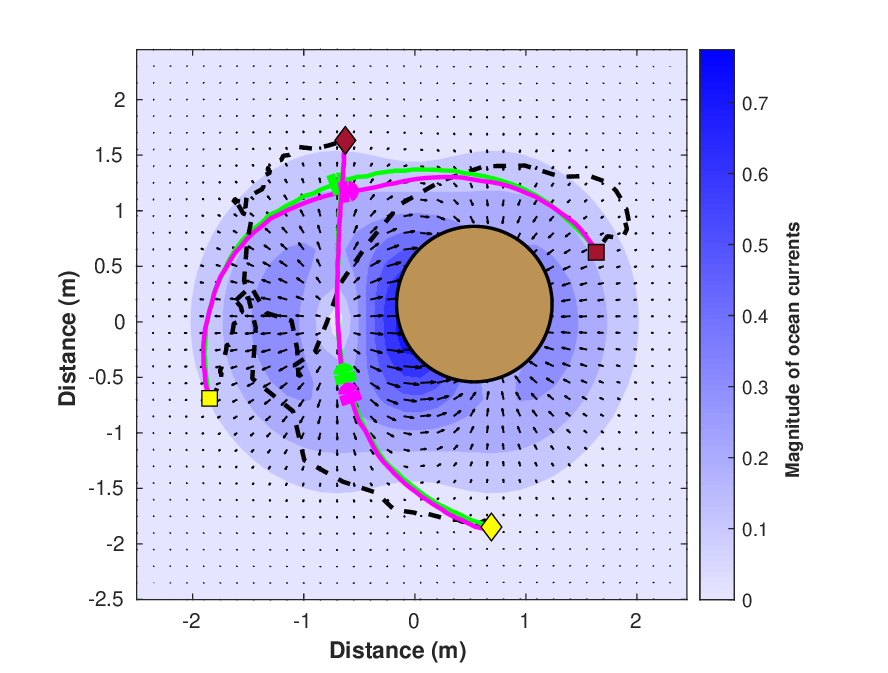}
\end{subfigure}
\hspace{-0.7cm}
\hfill
\hspace{-0.7cm}
\begin{subfigure}[b]{0.362\textwidth}
	\centering
	\includegraphics[width=\textwidth, trim={0.5cm 0 0 0},clip]{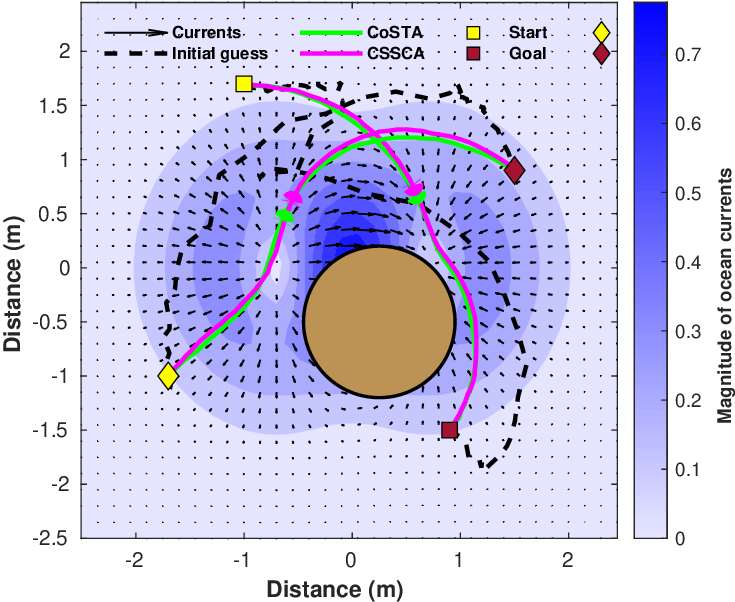}
\end{subfigure}
\newline
\vspace{-0.4cm}

\begin{subfigure}[b]{0.345\textwidth}
	\centering
	\includegraphics[width=\textwidth]{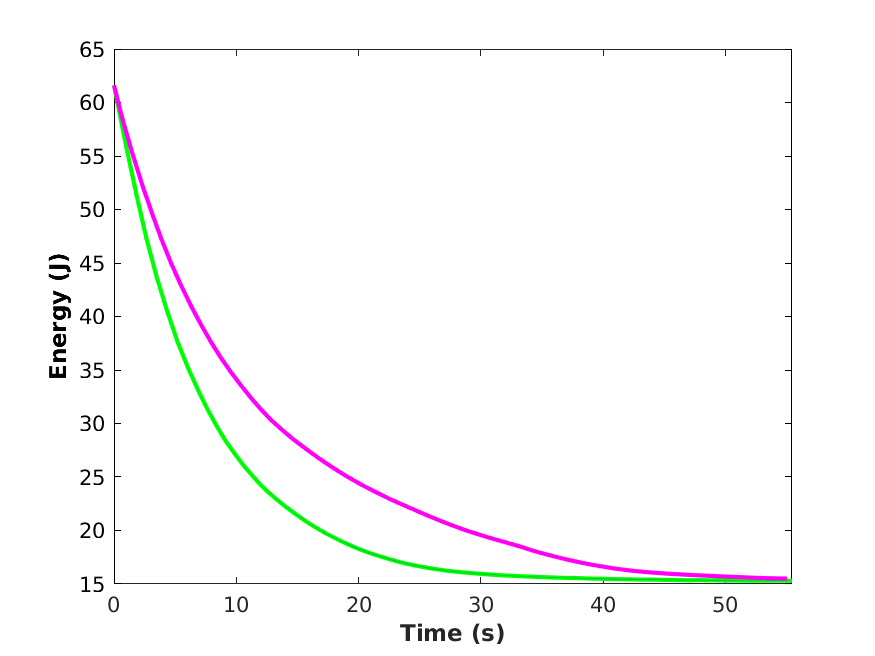}
	\caption{Environment 1}
\end{subfigure}
\hspace{-0.39cm}
\hfill
\hspace{-0.39cm}
\begin{subfigure}[b]{0.345\textwidth}
	\centering
	\includegraphics[width=\textwidth]{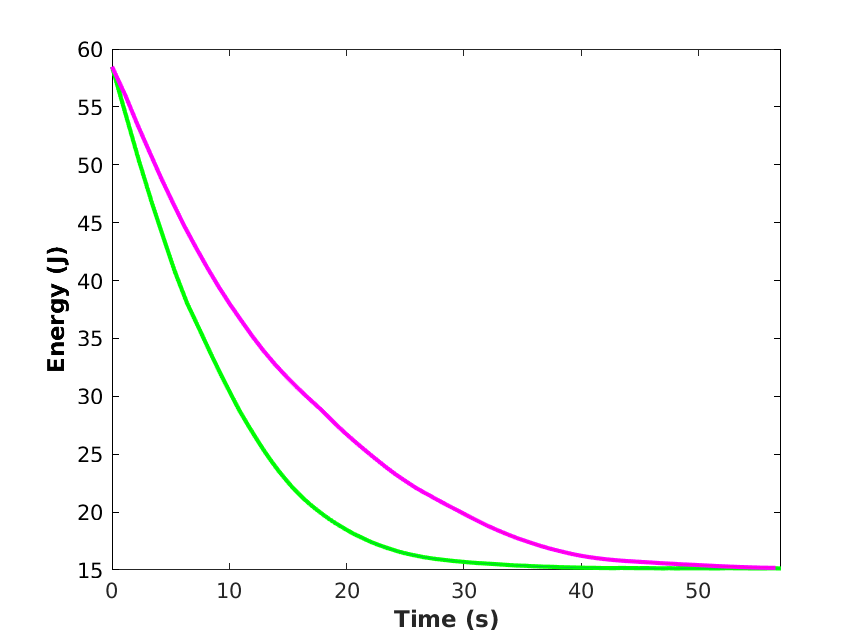}
	\caption{Environment 2}
\end{subfigure}
\hspace{-0.39cm}
\hfill
\hspace{-0.39cm}
\begin{subfigure}[b]{0.345\textwidth}
	\centering
	\includegraphics[width=\textwidth]{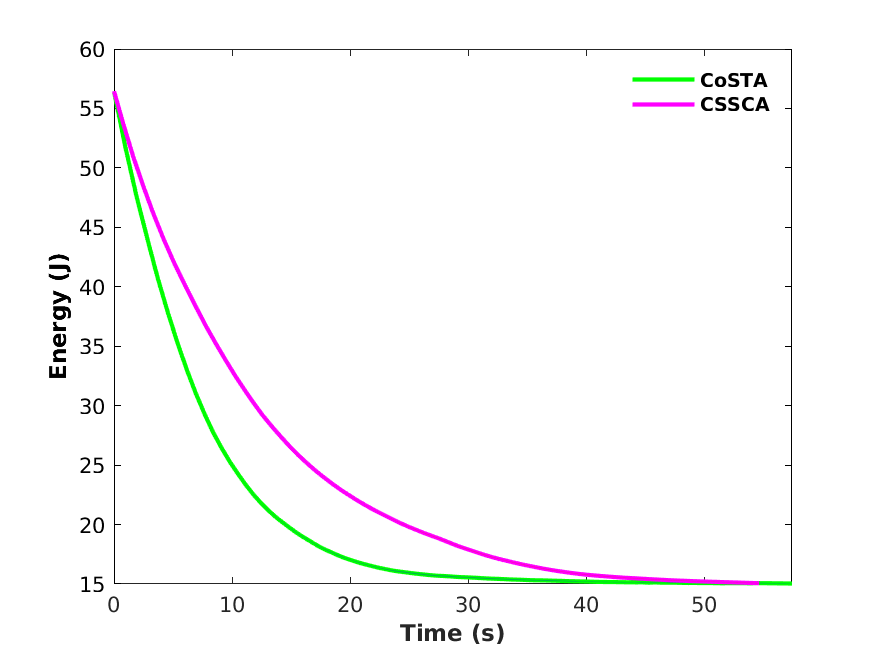}
	\caption{Environment 3}
	\label{figmp:f1c}
\end{subfigure}

\caption{Trajectory (top) and energy evolution vs time (bottom) plots for $T=60$, $T_f=15$, and $\sigma=0.2$, and environmental configurations 1, 2, and 3.}
\label{figmp:f1}
\end{figure*}

\begin{table}[]
\centering
	\caption{\label{ress}Classification accuracy of different algorithms for sparse logistic regression (where the entries are train accuracy/ test accuracy).}
	\scriptsize
	\begin{tabular}{c|c|c|c}
		\hline
		\textbf{Dataset}    & \textbf{Accuracy} & \textbf{LCPP} & \textbf{CoSTA} \\ 
		\hline
		\multirow[c]{2}{*}{MNIST}   & train    &  $90.5 \pm 1.4$    &   $94
		\pm 0.56 $         \\ \cline{2-4} 
		& test     & $90.54 \pm 1.51 $     &   $94.1
		\pm 0.74 $          \\ \hline
		\multirow[c]{2}{*}{gisette} & train    &  $99.7 \pm .05$    &   $99.44 \pm 0.11 $         \\ \cline{2-4} 
		& test     &  $97.36 \pm 0.5$    &  $97.72 \pm 0.32 $
		\\ 	\hline
	\end{tabular}
\end{table}

\begin{table}[]
	\centering
	\caption{Average time per iteration in seconds for $\sigma$ = 0.2, $T_f$ = 15, environmental configuration 4, and different values of $T$ for Case S1.}
	\label{tab:time}
	\begin{tabular}{lll}
		\hline
		& $T=30$ & $T=60$ \\ \cline{2-3} 
		\textbf{CoSTA} & 0.11$\pm$0.01 & 0.58$\pm$0.01 \\
		\textbf{CSSCA} & 0.10$\pm$0.01 & 0.59$\pm$0.01 \\ \hline
	\end{tabular}
\end{table}

\subsection{Energy-Optimal Trajectory Planning}\label{ssec:eematp}


Continuing the problem discussed in Sec. \ref{secc2}, we benchmark our proposed algorithm CoSTA against classical optimization based approach in CSSCA \cite{liu2019stochastic} and sampling based methods STOMP \cite{kalakrishnan2011stomp} and EESTO \cite{jones2017planning} on the energy-efficient trajectory planning problem \eqref{eqmp}-\eqref{eqmp:vr_max_s} in a simulated 2D ocean environment. Monte Carlo simulations were performed on an Intel(R) Xeon(R) E3-1226 CPU running Ubuntu 20.04 LTS with 32GB RAM. We modeled the environment and implemented the algorithms in MATLAB 2024a \cite{MATLAB23:online}.

We consider two cases of simulated ocean current setting with different complexity. For case \textit{S1} we consider a simple ocean current field with a single obstacle. The second case \textit{S2} considers navigation through a narrow passage with a complex ocean current field and multiple obstacles. For each case, the base ocean current field $\vartheta(\x)$ is generated differently. From the base ocean currents we generate current ensemble predictions $\vartheta(\x,\xib)$ by adding noise as follows: 
	\begin{align}\label{eqmp:nm}
		&\vartheta(\q,\xib) = \vartheta(\q)(\mathbf{I} + \text{diag}(\mathbf{e}(\xib)) ),
	\end{align}
	where $\mathbf{e(\xib)} \sim \mathcal{N} (\mathbf{0}, \sigma^2 \mathbf{I})$. As stated earlier, neither the true ocean currents nor the noise model is known to the vehicle, which can only query the predictor for noisy realizations $\vartheta(\q,\xib)$ of the currents. We discuss the setting and results for both the cases in detail in the following subsections.

\subsubsection*{Case S1 - Convergence Benchmark against optimization based algorithms}
	Here we consider a simple setting with one obstacle where the ocean currents at a point $\x = [x_1 ~ x_2]^\T$ are simulated using the following:
\begin{align} \label{oceancu}
	\vartheta(\x) = \omega \begin{bmatrix}
		1-2x_1^2 \\ 
		-2x_1x_2
	\end{bmatrix} \exp\{-(x_1^2+x_2^2)\}.
\end{align}

For S1 we compare the performance of proposed algorithm CoSTA with that of CSSCA across diverse settings. In  this case, we run simulations for multiple environments consisting of two agents and an obstacle. Each environment has a different starting, goal, and obstacle locations. Of these, the third environment, shown in Fig. \ref{figmp:f1c} is particularly challenging, as the ocean currents tend to push one of the agents into the obstacle, requiring more precise navigation. We also vary the values of number of way points $T$, mission time $T_f$, and the noise parameter $\sigma$.

The parameters of CSSCA ($\mu$, $\rho$,  and $\gamma$ \cite{liu2019stochastic}) and CoSTA ($\mu$, $\kb$, $c$, and $w$) were tuned for best possible performance after running the algorithms for initial 100 iterations. The results are averaged over 10 Monte Carlo runs for the following values of problem parameters: $\omega = 0.8$, $r^0 = 0.7$m, $r = 0.1$m, and $v^{\max} = 1$m/s. The initial trajectory for both the algorithms is generated by directly solving the non-convex feasibility problem \eqref{eqmp:term}-\eqref{eqmp:nc_obs2},\eqref{eqmp:vr_max_s} using the MATLAB's non-linear optimization function \texttt{fmincon}. The subproblems for CoSTA and CSSCA are Quadratic Constrained Quadratic Programs (QCQP) and solved using MATLAB function \texttt{fmincon} \cite{Linearor3:online} which uses the interior-point method to solve the convex subproblem. 

Fig. \ref{figmp:f1} shows the optimized trajectory of the agents along with evolution of energy cost vs CPU time\footnote{We omit the time required for first iteration for both algorithms from the calculations as we found it to be highly irregular. This is due to MATLAB Just In Time (JIT) compilation which may require more time in the first iteration due to generation and parsing of machine level code.} for various environments and configurations. For each of the three environments depicted, both algorithms ultimately generate similar trajectories. However, CoSTA always converges faster, achieving the (locally) optimal energy in less time. The faster convergence is attributed to momentum based gradient tracking in \eqref{stormtrack} and the adaptive step sizes in \eqref{etabeta} of CoSTA unlike that of CSSCA.

Table \ref{tab:time} shows the per iteration time for CoSTA and CSSCA for different number of waypoints $T$. As the number of waypoints increases, the size of each subproblem increases and hence more time is required. Note that since the subproblems for CoSTA and CSSCA are similar, their per-iteration time is also similar. Hence, we can conclude that the reason for speedup of CoSTA is that it requires fewer iterations due to its use of momentum. 

\subsubsection*{Case S2 - Effectiveness Benchmark against planning algorithms}
	For this case, we consider a more challenging scenario of narrow passage navigation through multiple obstaces. Currents are simulated using viscous lamb vortices \cite{sun2022efficient} which are generated using superposition of one point vortex solutions. The model with $M$ vortices is mathematically described by :
	\begin{equation}
		\vartheta(\x) = \sum_{m=1}^{M} \dfrac{\Omega_m(\x-\q_m)}{2\pi\norm{\x-\q_m}^2} \left( 1- \exp\left( -\dfrac{\norm{\x-\q_m}^2}{\delta_m^2}\right) \right)
	\end{equation}	
	with $\Omega_m = \begin{bmatrix}
		0 & -\omega_m \\ \omega_m & 0
	\end{bmatrix}$, where $\q_m$ is the centre of the $m$-th vortex and $\omega_m$, $\delta_m$ are parameters related to its strength and radius, respectively. In a 200m $\times$ 200m environment, we consider ocean model with $M=5$ Lamb vortices, $\omega \in \{-20,-10,20\}$, and $\delta \in \{10,20\}$. Two agents with $r = 5$m and $v^{\max} = 1$m/s seek to navigate to their respective goal locations in $T_f = 300s$ avoiding each other and 4 obstacles with radii $r^o \in \{10,20,30,35\}$ m.  

The initial trajectory for CoSTA for this case is generated in a similar way to case S1 by directly solving the non-convex feasibility problem \eqref{eqmp:term}-\eqref{eqmp:nc_obs2},\eqref{eqmp:vr_max_s} but this time with an energy-informed-AStar\footnote{The energy-informed-AStar is simple AStar search with energy as travel cost}\cite{kularatne2016time} initialization. This results in a overall feasible initial trajectory for CoSTA.

	\begin{table}[]
		\centering
		\caption{Comparison for Case S2}
		\label{tab:case-s2}
		\begin{tabular}{lcc}
			\hline
			Algorithm & \textbf{Energy (J)} & $\mathbb{E}$[\textbf{Constraint Violation}] \\ \hline 
			Initial guess & 243.71$\pm$1.53 & 0 \\
			CoSTA & \textbf{139.33$\pm$4.20} & \textbf{0} \\
			EESTO \cite{jones2017planning} & 604.84$\pm$22.02 & 6.94 \\
			STOMP \cite{kalakrishnan2011stomp} & 622$\pm$32.53 & 5.89 \\ \hline
		\end{tabular}
	\end{table}

\begin{figure}[]
	\centering
	\includegraphics[width=0.48\textwidth]{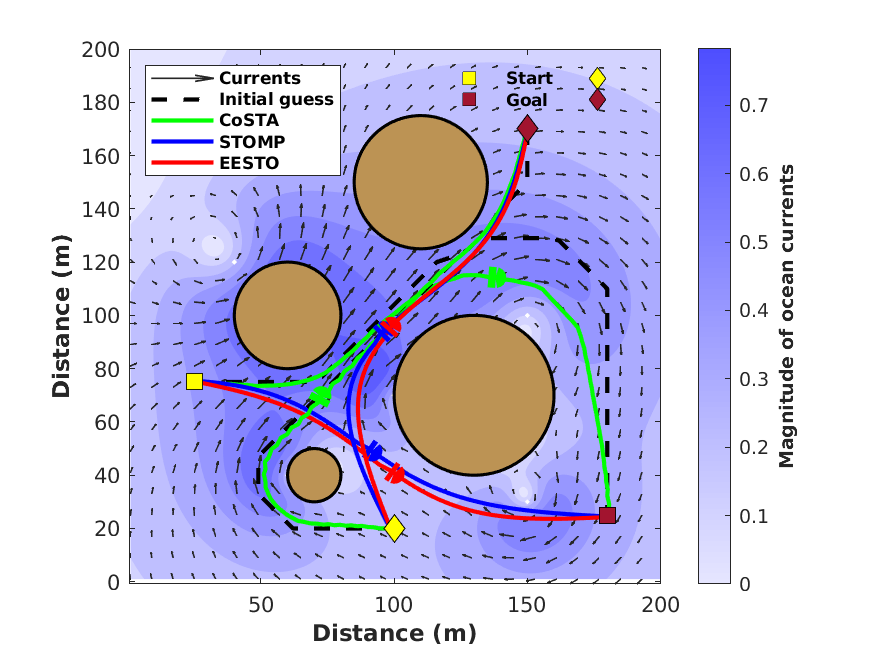}
	
	\caption{Trajectory plots for case S2 environment with 5 Lamb vortices.}\label{figmp:lamb-a}
\end{figure}

We benchmark CoSTA against sampling-based planning methods STOMP and EESTO. All the algorithms were independently tuned for best possible performance. For this case all algorithms are again run for 10 Monte Carlo runs. A sample trajectory obtained from one of the runs is shown in Fig. \ref{figmp:lamb-a}. Observe, in Fig. \ref{figmp:lamb-a}, for one of the agents CoSTA trajectory takes a long path along the flow to conserve energy whereas STOMP and EESTO get stuck in local minimum taking a shorter but flow opposing energy hungry trajectory. Sampling-based approaches such as STOMP and EESTO operate using average ocean current fields, as they are not equipped to incorporate ensemble current forecasts. Consequently, their performance tends to be suboptimal in this setting. 

Table \ref{tab:case-s2} reports the average energy consumption and average constraint violations over ocean ensembles for the competing algorithms. All algorithms produce collision-free trajectories; however, CoSTA achieves the lowest average energy consumption without any constraint violations, thereby demonstrating superior performance.

\section{Conclusion} \label{conc}
We have considered the problem of minimizing a stochastic non-convex optimization problem with non-convex constraints. While the problem has been well-studied, the state-of-the-art algorithms have suboptimal oracle complexity of $\O(\epsilon^{-2})$. The recently proposed algorithm in \cite{shi2022momentum} achieves the optimal rate of $\O(\epsilon^{-3/2})$, but this result is only valid under the condition that the initial point is nearly feasible  with non-adaptive step sizes. In contrast, we propose the CoSTA algorithm, which introduces recursive momentum updates to achieve  an almost optimal oracle complexity of $\Ot(\epsilon^{-3/2})$ under adaptive step sizes and Lipschitz smoothness assumptions. Furthermore, under a milder Mean Squared Smoothness (MSS) condition, CoSTA retains this complexity with non-adaptive step sizes. The proposed algorithm is developed within the SCA framework, and requires convex surrogates of the objective and constraint functions at every iteration. We also introduced a parameterized version of the Mangasarian-Fromovitz constraint qualification that allows us to ensure that the intermediate subproblems are strictly feasible and the intermediate dual variable stay bounded, resulting in an oracle complexity bound that depends only on the problem parameters. Finally, the empirical performance of the  CoSTA algorithm is compared with the other state-of-the-art algorithms across two different applications: energy optimal trajectory design and sparse binary classification. The proposed algorithm outperforms the state-of-the-art algorithms on the energy optimal trajectory design problem and performs at par with the LCPP algorithm on the sparse binary classification problem. 

\appendices

\section{Proof of Lemma \ref{slater}} \label{ap3}
\begin{IEEEproof}
	Since $\x_t \in \cX$, parameterized MFCQ holds at $\x_t$, implying the existence of $\d_t$ as specified in Assumption \eqref{Amfcq}. Consider the point $\tx(\x_t) = \x_t + \gamma \d_t$ for $\gamma > 0$ and let $\cI$ and $\cJ$ be defined as the sets of indices corresponding to the near-active constraints at $\x_t$ for the convex and general constraints, respectively. 
	
	First consider the near-active constraints at the feasible point $\x_t$. From the smoothness of $\gt_j(\cdot,\x_t)$, we have that
	\begin{align}
		\gt_j(\tx(\x_t),\x_t) &\leq \gt_j(\x_t,\x_t) + \ip{\nabla \gt_j(\x_t,\x_t)}{\tx(\x_t)-\x_t} + \frac{L}{2}\norm{\tx(\x_t)-\x_t}^2 \\
		&\hspace{-1.5cm}=\gt_j(\x_t) + \gamma\ip{\nabla \gt_j(\x_t)}{\d_t} + \frac{\gamma^2L}{2} \leq - \rho\gamma + \frac{\gamma^2L}{2} 
	\end{align}
	for all $j \in \cJ$.  Likewise, from the smoothness of $h_i$, we have that 
	\begin{align}
		h_i(\tx(\x_t))  &\leq - \rho\gamma + \frac{\gamma^2L}{2} 
	\end{align}
	for all $i \in \cI$. Therefore, if we set $\gamma = \frac{\rho}{L}$, it can be seen that $\gt_j(\tx(\x_t),\x_t) \leq -\frac{\rho^2}{2L}$ and $h_i(\tx(\x_t)) \leq -\frac{\rho^2}{2L}$ as required. Since $\tx(\x_t) \in \cX$, it follows from Assumption \ref{cset} that $\tx(\x_t) \in \cD$ and both the constraint functions are well-defined at $\tx(\x_t)$. 
	
	Next, let us consider the strongly inactive constraints, where we have that
	\begin{align}
		\gt_j(\tx(\x_t),\x_t) &= \gt_j(\tx(\x_t),\x_t) - \gt_j(\x_t,\x_t) + g_j(\x_t) \\
		&\leq G\norm{\tx(\x_t)-\x_t} - \chi= \frac{G\rho}{L} - \chi 
	\end{align} 
	for all $j \notin \cJ$, and likewise, from the Lipschitz continuity of $h_i$, we have that 
	\begin{align}
		&h_i(\tx(\x_t)) \leq \frac{G\rho}{L} - \chi 
	\end{align}  
	for all $i \notin \cI$. Therefore, the required condition is satisfied if we have $\chi \geq \frac{\rho}{L}\left(G + \rho/2\right)$. 
\end{IEEEproof}

\section{Proof of Lemma \ref{dualvari}} \label{ap4}
\begin{IEEEproof}
	Since $\ft(\cdot,\x_t,\z_t,\xib_t)$, $\gt(\cdot,\x_t)$, and $h(\cdot)$ are convex functions, we have that for any $\x \in \cX$, 
	\begin{align}
		\ft(\x,\x_t,\z_t,\xib_t)&+u(\x) \geq \ft(\xhat,\x_t,\z_t,\xib_t) + u(\xhat)\nonumber  \\
		&+ \ip{\nabla \ft(\xhat,\x_t,\z_t,\xib_t)+\v_{\xhat}}{\x-\xhat} \label{ftconv}\\
		\gt_j(\x,\x_t) &\geq \gt_j(\xhat,\x_t) + \ip{\nabla \gt_j(\xhat,\x_t)}{\x-\xhat} \label{gjconv} \\
		h_i(\x) &\geq h_i(\xhat) + \ip{\nabla h_i(\xhat)}{\x-\xhat} \label{hiconv} 
	\end{align}
	for all $1\leq j \leq J$ and $1\leq i \leq I$. Setting $\x = \tx(\x_t)$, multiplying $\eqref{gjconv}$ by $[\lamh_t]_j \geq 0$, multiplying \eqref{hiconv} by $[\hat{\nu}_t]_i$, and adding with \eqref{ftconv}, we obtain
	\begin{align}
		&\ft(\tx(\x_t),\x_t,\z_t,\xib_t) + u(\tx(\x_t)) \nonumber \\
		& + \sum\nolimits_{j=1}^J [\lamh_t]_j\gt_j(\tx(\x_t),\x_t) + \sum\nolimits_{i=1}^I [\hat{\nu}_t]_i h_i(\tx(\x_t)) \nonumber\\
		&\geq \ft(\xhat,\x_t,\z_t,\xib_t) + u(\xhat) \nonumber \\
		&~~~+ \sum\nolimits_{j=1}^J [\lamh_t]_j\gt_j(\xhat,\x_t)  + \sum\nolimits_{i=1}^I [\hat{\nu}_t]_i h_i(\xhat) \nonumber \\
		&+ \Big \langle\sum\nolimits_{j=1}^J[\lamh_t]_j\nabla\gt_j(\xhat,\x_t) + \sum\nolimits_{i=1}^I [\hat{\nu}_t]_i\nabla h_i(\xhat) \nonumber \\
		&~~~+  \nabla \ft(\xhat,\x_t,\z_t,\xib_t) + \v_{\xhat},\tx(\x_t)  - \xhat\Big \rangle  \label{prf22}
	\end{align}
	The last three terms on the right vanish from \eqref{cs} and \eqref{stat2}. So using Lemma \ref{slater}, we have that
	\begin{align}
		&\ft(\tx(\x_t),\x_t,\z_t,\xib_t) + u(\tx(\x_t)) - \ft(\xhat,\x_t,\z_t,\xib_t) - u(\xhat) \nonumber\\
		&\geq - \sum\nolimits_{j=1}^J [\lamh_t]_j\gt_j(\tx(\x_t),\x_t) -\sum\nolimits_{i=1}^I [\hat{\nu}_t]_i h_i(\tx(\x_t)) \label{dualmid}\\
		&\geq  \frac{\rho^2}{2L}\bigg(\sum_{j=1}^J [\lamh_t]_j+\sum_{i=1}^I [\hat{\nu}_t]_i \bigg). \label{prf33}
	\end{align}
	Finally, using Assumption \eqref{Asul2} to upper bound the terms on the left, we obtain the desired result. 
\end{IEEEproof}

	\section{Proof of Lemma \ref{stproof}} \label{ap1} 
\begin{IEEEproof}	
	The proof follows along the lines of \cite[Lemma 2]{cutkosky2019momentum}. Expanding $\z_{t+1} - \nabla U(\x_t)$ from \eqref{stormtrack} and introducing $ (1-\beta_{t+1})\nabla U(\x_{t})$, we obtain
	\begin{align}
		\z_{t+2} -& \nabla U(\x_{t+1}) = (1-\beta_{t+1})(\z_{t+1}  - \nabla U(\x_t)) \nonumber\\
		& - (1-\beta_{t+1})(\nabla f(\x_t,\xib_{t+1})-\nabla U(\x_t) ) \nonumber\\
		&+ \nabla f(\x_{t+1},\xib_{t+1}) - \nabla U(\x_{t+1}) 
	\end{align}
	where observe that $\EE_{\xib_{t+1}}[\nabla U(\x_t) - \nabla f(\x_t,\xib_{t+1})] = \EE_{\xib_{t+1}}[\nabla f(\x_{t+1},\xib_{t+1}) - \nabla U(\x_{t+1})] = 0$. Therefore, we have that
	\begin{align}
		&\EE_{\xib_{t+1}}[\tfrac{\norm{\z_{t+2} - \nabla U(\x_{t+1})}^2}{\eta_t}] \nonumber \\
		&~ = \EE_{\xib_{t+1}} \bigg[\frac{(1-\beta_{t+1})^2\norm{\z_{t+1} - \nabla U(\x_t)}^2}{\eta_t} \bigg ] \nonumber\\
		&~~~~ +\EE_{\xib_{t+1}} \Big[ \eta_t^{-1}\big\| (1-\beta_{t+1}) \big (\nabla f(\x_t,\xib_{t+1}) \nonumber \\
		&~~~~~~~~-\nabla U(\x_t)+ \nabla U(\x_{t+1}) - \nabla f(\x_{t+1},\xib_{t+1}) \big )\nonumber \\
		&~~~~~~~~~~~ - \beta_{t+1}\big ( \nabla f(\x_{t+1},\xib_{t+1}) - \nabla U(\x_{t+1})) \big\|^2\Big]  \label{stproof2}
	\end{align}
	where the cross term vanishes since $\z_{t+1} - \nabla U(\x_t)$ is independent of $\xib_{t+1}$ and the second summand is zero mean. The second term in \eqref{stproof2} can again be expanded as
	\begin{align}
		&\EE_{\xib_{t+1}} \Big[ \eta_t^{-1}\big\| (1-\beta_{t+1}) \big (\nabla f(\x_t,\xib_{t+1}) \nonumber \\
		&~~~~~-\nabla U(\x_t)+ \nabla U(\x_{t+1}) - \nabla f(\x_{t+1},\xib_{t+1})\big )\nonumber \\
		&~~~~~~~~~~~ - \beta_{t+1} \big ( \nabla f(\x_{t+1},\xib_{t+1}) - \nabla U(\x_{t+1}) ) \big\|^2\Big]\nonumber \\
		& \overset{(a)}{\leq} 2\EE_{\xib_{t+1}} \bigg[\frac{(1-\beta_{t+1})^2\norm{\nabla f(\x_{t+1},\xib_{t+1}) - \nabla f(\x_t,\xib_{t+1})}^2}{\eta_t}\bigg] \nonumber \\
		&~~~~~~~~~~ + 2\EE_{\xib_{t+1}} \bigg[\frac{\beta_{t+1}^2\norm{\nabla f(\x_{t+1},\xib_{t+1}) - \nabla U(\x_{t+1})}^2}{\eta_t}\bigg] \nonumber \\
		& \overset{(b)}{\leq} 2\EE_{\xib_{t+1}} \bigg[\frac{(1-\beta_{t+1})^2\norm{\nabla f(\x_{t+1},\xib_{t+1}) - \nabla f(\x_t,\xib_{t+1})}^2}{\eta_t}\bigg] \nonumber \\
		&~~~~~~~~~~ + 2\EE_{\xib_{t+1}} \bigg[\frac{\beta_{t+1}^2\norm{\nabla f(\x_{t+1},\xib_{t+1}) }^2}{\eta_t}\bigg]  \label{stproof3}
	\end{align}
	where in (a) we have used the inequality $\EE[\norm{\mathsf{X}-\EE[\mathsf{X}] + \mathsf{Y}}^2] \leq 2\EE[\norm{\mathsf{X}}^2] + 2\EE[\norm{\mathsf{Y}}^2]$ for any random variables $\mathsf{X}$ and $\mathsf{Y}$ with bounded variances. And in (b) we have used $\EE[\norm{\mathsf{X}-\E{\mathsf{X}}}^2] \leq \EE[\norm{\mathsf{X}}^2]$. 	Dividing both sides by $\tfrac{1}{16 L^2}$, taking full expectation in  \eqref{stproof2} and using $\beta_{t+1} \leq \frac{1}{4}$, \eqref{stproof2}  can be rewritten as, 		
	%
	%
	\begin{align}
		&\frac{1}{16L^2}  \E{\frac{\norm{\e_{t+1}}^2}{\eta_t}} \leq \frac{1}{16L^2} \E{\frac{(1-\beta_{t+1})\norm{\e_t}^2}{\eta_t}} \nonumber \\
		&~~+ \frac{1}{8L^2} \E{\frac{\norm{\nabla f(\x_{t+1},\xib_{t+1}) - \nabla f(\x_t,\xib_{t+1})}^2}{\eta_t}} \label{inte1} \\
		&~~~~ + \frac{1}{8L^2}\E{\frac{\beta_{t+1}^2G_{t+1}^2}{\eta_t}} \label{inte}
	\end{align}
	Using $L$-smoothness, and the update equation \eqref{stormupdate}, we obtain
	\begin{align}
		\frac{1}{16L^2} \E{\frac{\norm{\e_{t+1}}^2}{\eta_t}} 
		& \leq  \frac{1}{16L^2} \E{\frac{(1-\beta_{t+1})\norm{\e_t}^2}{\eta_t}} + \frac{1}{8}\E{\frac{\norm{\x_{t+1} - \x_t}^2}{\eta_t}} \nonumber \\
		&~~~~~~~~~ + \frac{1}{8L^2}\E{\frac{\beta_{t+1}^2 G_{t+1}^2}{\eta_t}}   \\
		&\leq \frac{1}{16L^2} \E{\frac{(1-\beta_{t+1})\norm{\e_t}^2}{\eta_t}} + \frac{1}{8}\E{\eta_t\norm{\del_t}^2} \nonumber \\
		&~~~~~~~~~+ \frac{1}{8L^2}\E{\frac{\beta_{t+1}^2 G_{t+1}^2}{\eta_t}}. \label{stproof5}
	\end{align}

\end{IEEEproof}

\section{Proof of Lemma \ref{le4}} \label{ap2}
\begin{IEEEproof}
	We first obtain a fundamental inequality. The optimality condition of \eqref{xhata} implies that
	\begin{align}
		&\big \langle \y - \xhat,\nabla \ft(\xhat,\x_t,\z_t,\xib_t) + \v_\xhat \big \rangle \geq  0
	\end{align}
	for all $\y \in \cX(\x_t)$, where $\v_\xhat \in \partial u(\xhat)$. Choosing $\y = \x_t$,
	\begin{align}
		&\big \langle \x_t - \xhat,\nabla \ft(\xhat,\x_t,\z_t,\xib_t)+ \v_\xhat \big \rangle \geq  0.
	\end{align}
	Adding and subtracting $\nabla \ft(\x_t,\x_t,\z_t,\xib_t) = \z_{t+1}$, we obtain
	\begin{align}
		&\ip{\xhat - \x_t}{\z_{t+1} + \v_\xhat} \label{mida}\\
		&~ + \langle \xhat -\x_t,\nabla \ft(\xhat,\x_t,\z_t, \xib_t) - \nabla \ft(\x_t,\x_t,\z_t,\xib_t) \rangle \leq 0 \nonumber, \\
		&\Rightarrow \ip{\z_{t+1} + \v_\xhat}{\xhat -\x_t} + \mu \norm{\xhat - \x_t}^2 \leq 0, \label{midb}  
	\end{align}
	where we have used the strong monotonicity property of the gradient of the strongly convex function $\ft(\cdot,\x_t,\z_t,\xib_t)$. 
	%
	From the convexity of $u$ and the update rule \eqref{stormupdate}, we obtain
	\begin{align}
		u{(\x_{t+1})} -  u(\x_t)  &\leq \eta_t \ip{\v_\xhat}{\xhat - \x_t} 
	\end{align}
	Multiplying \eqref{midb} by $\eta_t$ and adding, we obtain
	\begin{align}
		&\eta_t\ip{\z_{t+1}}{\xhat -\x_t} \leq - \mu\eta_t \norm{\xhat - \x_t}^2 + u(\x_t) - u(\x_{t+1}).
	\end{align}
	Adding and subtracting $\nabla U(\x_t)$, we obtain
	\begin{align}
		&\eta_t\ip{\z_{t+1} -\nabla U(\x_t)}{\xhat -\x_t} \nonumber \\
		&~~~ +\eta_t\ip{\nabla U(\x_t)}{\xhat -\x_t}  \nonumber\\
		&~~~\leq - \mu \eta_t \norm{\xhat - \x_t}^2   + u(\x_t) - u(\x_{t+1}). \label{ht}
	\end{align}		
	We now proceed to bound $U(\x_{t+1})$ in terms of $U(\x_t)$. Since $U$ is smooth, we have that:
	\begin{align}
		&U(\x_{t+1}) \leq U(\x_t) + \ip{\nabla U(\x_t)}{ \x_{t+1}-\x_t} + \frac{L}{2}\norm{\x_{t+1}-\x_t}^2\\
		&\leq U(\x_t) + \eta_t \ip{\nabla U(\x_t)}{\xhat - \x_t}+ \frac{L\eta_t^2}{2}\norm{\xhat - \x_t}^2\\
		&\leq U(\x_t) + u(\x_t) - u(\x_{t+1}) - \eta_t\ip{\z_{t+1} -\nabla U(\x_t)}{\xhat -\x_t} \nonumber\\
		&~~~  + \bigg(\frac{L}{2}\eta_t^2- \mu \eta_t\bigg)\norm{\xhat -\x_t}^2
	\end{align}
	where the last inequality follows from substituting \eqref{ht}. Using  Young's inequality, rearranging, and taking expectation on both sides, we obtain
	\begin{align}
		&\E{F(\x_{t+1}) - F(\x_t)} \leq  \frac{1}{2}\E{\eta_t\norm{\nabla U(\x_t) - \z_{t+1} }^2} \nonumber \\
		&~~~+ \frac{L}{2}\E{\eta_t^2\norm{\xhat - \x_t}^2} +\bigg(   \frac{1}{2} - \mu \bigg) \E{\eta_t\norm{\xhat - \x_t}^2}  \nonumber \\
		&\leq\frac{1}{2}\E{\eta_t \norm{\e_t}^2} - \frac{1}{4}\E{(4\mu  - 2L\eta_t - 2  ) \eta_t \norm{\del_t}^2} \label{stdescent}
	\end{align}
	Finally, using $\mu >\frac{ L\eta_t}{2} + \frac{3}{4}$, we get the desired result.
\end{IEEEproof}

\section{Proof of Corollary \ref{cor2} and Corollary \ref{cor3}} \label{prc2c3}
\begin{IEEEproof}
	Under non-adaptive step sizes, specifically $\eta_t = \tfrac{\kb}{(w+t)^{1/3}}$ and $\beta_{t+1} = c\eta_t^2$. The results of Lemma \ref{stproof} and \ref{le4} the step sizes $\eta_t$ can be taken outside the expectation. Consequently, the modified result in Lemma \ref{stproof} becomes,		
	\begin{align} \label{rs1}
		\frac{1}{16L^2 \eta_t} &\E{\norm{\e_{t+1}}^2} \leq \frac{1 - \beta_{t+1}}{16 L^2 \eta_t} \E{\norm{\e_t}^2} \nonumber \\
		&+ \frac{\eta_t}{8} \E{\norm{ \del_t}^2} + \frac{\beta_t^2}{8L^2 \eta_t} \E{(G_{t+1})^2}.
	\end{align}
	Lemma \ref{le4} is also modified as follows.		
	\begin{align} \label{rs2}
		\E{F(\x_{t+1}) - F(\x_t)} \leq \frac{\eta_t}{2} \E{\norm{\e_t}^2} - \frac{\eta_t}{4} \E{\norm{\del_t}^2}.
	\end{align}	
	Note that here only Assumption \ref{mss} is required instead of Assumption \ref{assmooth} (as $\eta_t $ can taken outside the expectation in \eqref{inte1}). Following the proof structure of Theorem \ref{th1}, and introducing the Lyapunov function $\Psi_t = \E{F(\x_t) - F^\star} + \tfrac{1}{16 L^2 \eta_{t-1}} \E{\norm{\e_t}^2}$, we obtain,
	\begin{align} \label{stproof11}
		\Psi_{t+1} - \Psi_t &\leq \frac{\alpha_t}{16 L^2} \E{\norm{\e_t}^2} - \frac{\eta_t}{8} \E{\norm{\del_t}^2} + \frac{c^2 \eta_t^3}{8L^2} \E{G_{t+1}^2},
	\end{align}
	where $\alpha_t := 8L^2 \eta_t + \frac{1 - \beta_{t+1}}{\eta_t} - \frac{1}{\eta_{t-1}}$. Just like \eqref{th1proof1}, we can again write 
	\begin{align}
		\frac{\kb}{\eta_t} - \frac{\kb}{\eta_{t-1}}& =(w +t)^{1/3} - (w+ t-1)^{1/3}  \nonumber\\
		&\leq \frac{1}{3(w+t-1)^{2/3}} = \frac{\eta_{t-1}^2}{3\kb^2}.
	\end{align}
	Now, similar to \eqref{th1proof2}-\eqref{th1proof3}, we have
	\begin{align}
		\alpha_t &\coloneq 8 L^2\eta_t +	\frac{(1-\beta_{t+1})}{\eta_t} - \frac{1}{\eta_{t-1}} \nonumber \\
		&= 8 L^2\eta_t  + \frac{1}{\eta_t} - \frac{1}{\eta_{t-1}} - \frac{\beta_{t+1}}{\eta_t} \\
		&\leq 8 L^2\eta_{t-1}  + \frac{\eta_{t-1}^2}{3\bar{k}^3} - c\eta_{t-1} \leq -2L^2 d' \eta_{t-1}\label{th1proof41}
	\end{align}
	where $d':= \tfrac{1}{2L^2}(c - 8L^2 - \tfrac{1}{3\bar{k}^3}) > 0$. Therefore the first term on the right of \eqref{stproof11} can be dropped. Taking telescopic sum of \eqref{stproof11} from $t = 1 \ldots T$ and rearranging,
	\begin{align}
		\frac{1}{8} \sum_{t =1}^T \eta_t\E{\norm{\del_t}^2 }  &\leq \Psi_1 - \Psi_{T+1} + \frac{c^2 G^2  }{8L^2  } \sum_{t =1}^T \eta_t^3  \label{rearg}
	\end{align}
	Similar to  \eqref{logsum} we can bound, 
	\begin{align}
		& \sum\nolimits_{t=1}^T\eta_t^3   =  \kb^3   \sum\nolimits_{t=1}^T \frac{1}{w+ t}  \leq  \kb^3 \ln(1 + \frac{T}{w}) \nonumber \\
		&\qquad=   \mathcal{O}(\ln(1 + \frac{T}{w})) 
	\end{align}
	Like \eqref{psibound} we can bound,
	\begin{align}
		\Psi_1 \leq F(\x_1) - F^\star + \frac{\sigma^2}{16L^2\eta_0}= B_1 + \frac{w^{1/3}\sigma^2}{16L^2 \kb} . 
	\end{align}
	Hence, finally, using similar inequalities as  \eqref{delta} we lower bound \eqref{rearg},
	\begin{subequations}
		\begin{align}
			\Delta_T^2 \leq M'_T\frac{(w+T)^{1/3}}{T} 
		\end{align}
		where $M'_T := 8\frac{B_1}{\kb} + \frac{w^{1/3}\sigma^2}{2L^2\bar{k}^2} + \frac{ c^2\kb^2 G^2}{ L^2}  \ln (1 + \frac{T}{w})$. 
		Similarly, we get
		
		\begin{align}
			\tfrac{1}{T} \sum_{t =1}^{T}  \E{\norm{\del_t}^2}  \leq \frac{M'_T(w+T)^{1/3}}{T} = \Ot \left(\frac{1}{T^{2/3}}\right) \label{cor:eq1} \\
			\tfrac{1}{T}\sum_{t=1}^T\E{\norm{\e_t}} \leq \sqrt{\frac{M'_T(w+T)^{1/3}}{d'T}} = \Ot \left(\frac{1}{T^{1/3}}\right) \label{cor:eq2}
		\end{align}
	\end{subequations}
	where 	$d':= \tfrac{1}{2L^2}(c - 8L^2 - \tfrac{1}{3\bar{k}^3}) > 0$. 	Hence, we obtain the desired result in Corollary \ref{cor2}. 
	
	Also, under non-adaptive diminishing step size, in light of Corollary \ref{cor2}, \eqref{frate} can be bounded as,
	\begin{subequations}
		\begin{align}
			&\tfrac{1}{T}\sum\nolimits_{t=1}^T \E{\norm{\pib_t}}^2 \nonumber \\
			&~~~~ \leq \tfrac{2}{T}\sum_{t=1}^T \left(2L + \frac{4B_UL^2}{\rho^2}\right)^2  \E{\norm{\del_t}}^2+ \tfrac{2}{T}\sum_{t=1}^T \E{\norm{\ep_t}}^2 \nonumber \\
			&~~~~\leq \Ot (T^{-2/3})  \\
			&\text{Similarly, in light of Corollary \ref{cor2},  \eqref{cs1} can bounded as,} \nonumber \\
			& \tfrac{1}{T}\sum_{t =1}^{T} \EE [\lamh_t^\T\g(\xhat)] \geq -\frac{2B_U L^2}{\rho^2 T} 	\sum_{t =1}^{T}   \EE [\norm{\del_t}^2]  = \Ot (T^{-2/3}).  
		\end{align}
		Therefore,
		\begin{align}
			&\min_{1\leq t\leq T} (\E{\norm{\pib_t}}^2 - \EE [\lamh_t^\T\g(\xhat)]) \nonumber \\
			&~~~~~~~~~\leq \tfrac{1}{T}\sum\nolimits_{t=1}^T \E{\norm{\pib_t}}^2  -   \tfrac{1}{T}\sum\nolimits_{t =1}^{T} \EE [\lamh_t^\T\g(\xhat)] \nonumber \\
			&~~~~~~~~~ \leq \Ot(T^{-2/3})  
		\end{align}
	\end{subequations} 
	Hence, we obtain the desired result in Corollary \ref{cor3}.
\end{IEEEproof}

\footnotesize
\bibliographystyle{ieeetr} 
\bibliography{IEEEabrv,cit}

\end{document}